\documentclass[aos,MSNbibl,seceqn,dvips]{arximspdf}
\usepackage{mathbh,dcolumn}
\usepackage{graphicx}

\doi{10.1214/12-AOS1038} 
\volume{40}
\issue{4}
\pubyear{2012}
\firstpage{2327}
\lastpage{2355}

\makeatletter
\newcolumntype{d}[1]{D{.}{.}{#1}}
\newcommand{\eqref}[1]{(\ref{#1})}
\newtheorem{theorem}{Theorem}

\newproclaim{definition}{Definition}
\newtheorem{lemma}{Lemma}
\newtheorem{proposition}{Proposition}
\newproclaim{remark}{Remark}

\newcommand{\R}{\mathbb{R}}
\newcommand{\N}{\mathbb{N}}
\newcommand{\E}{\mathbb{E}}

\newcommand{\nnabla}{\bolds{\nabla}}
\newcommand{\YY}{\mathbf{Y}}
\newcommand{\bZ}{\mathbf{Z}}
\newcommand{\bb}{\mathbf{b}}
\newcommand{\ff}{\mathbf{f}}
\newcommand{\boldg}{\mathbf{g}}
\newcommand{\ssigma}{\bolds{\sigma}}
\newcommand{\KL}{\mathcal{K}}
\def\xxi{\bolds{\xi}}
\def\bbeta{\bolds{\beta}}
\def\llambda{\bolds{\lambda}}
\def\ttheta{\bolds{\theta}}
\def\bu{\mathbf{u}}
\def\bv{\mathbf{v}}

\newcommand{\sgn}{\operatorname{sgn}}
\newcommand{\diag}{\operatorname{diag}}
\newcommand{\Oracle}{\operatorname{Oracle}}
\newcommand{\Tr}{\operatorname{Tr}}
\newcommand{\Range}{\operatorname{Range}}
\def\DST{\mathcal{D}}
\newcommand{\MSE}{\operatorname{MSE}}

\newcommand{\1}{\mathbh{1}}

\newcommand{\eg}{e.g., }
\newcommand{\lcf}{cf. }

\makeatother

\begin{document}
\begin{frontmatter}

\title{Sharp oracle inequalities for aggregation of affine~estimators\thanksref{T1}}
\runtitle{Aggregation of affine estimators}

\begin{aug}
\author[a]{\fnms{Arnak S.}~\snm{Dalalyan}\corref{}\ead[label=e1]{arnak.dalalyan@ensae.fr}\ead[label=u1,url]{http://www.arnak-dalalyan.fr/}}
\and
\author[b]{\fnms{Joseph}~\snm{Salmon}\ead[label=e2]{salmon@math.jussieu.fr}\ead[label=u2,url]{http://josephsalmon.eu/}}
\thankstext{T1}{Supported in part by ANR Parcimonie.}
\runauthor{A.~S. Dalalyan and J.~Salmon}
\affiliation{ENSAE-Crest, Universit\'e Paris Est and Universit\'e Paris Diderot}
\address[a]{ENSAE-Crest\\
3 Avenue Pierre Larousse\\
92245 Malakoff Cedex\\
France
\printead{e1}\\
\printead{u1}}

\address[b]{2 Pl. Jussieu BP 7012\\
75251 PARIS Cedex 05\\
France \\
\printead{e2}\\
\printead{u2}}
\end{aug}

\received{\smonth{4} \syear{2011}}
\revised{\smonth{6} \syear{2012}}

%
\begin{abstract}
We consider the problem of combining a (possibly uncountably infinite)
set of affine estimators in
nonparametric regression model with heteroscedastic Gaussian
noise. Focusing on the exponentially weighted aggregate, we prove a
PAC-Bayesian type inequality that leads
to sharp oracle inequalities in discrete but also in
continuous settings. The framework is general enough to cover the
combinations of various procedures such as
least square regression, kernel ridge regression,
shrinking estimators and many other estimators used in the literature
on statistical inverse problems. As a
consequence, we show that the proposed aggregate
provides an adaptive estimator in the exact minimax sense without
discretizing the range of tuning
parameters or splitting the set of observations. We also illustrate
numerically the good performance achieved
by the exponentially weighted aggregate.
\end{abstract}

%


\end{frontmatter}
%
\section{Introduction}
There is growing empirical evidence of superiority of aggregated
statistical procedures, also referred to as
\textit{blending},
\textit{stacked generalization} or \textit{ensemble methods}, with
respect to ``pure'' ones. Since their
introduction in the 1990s, famous aggregation procedures such as
\textit{Boosting}~\cite{Freund90}, \textit{Bagging}~\cite{Breiman96b}
or \textit{Random Forest}~\cite{AmitGeman97} have been successfully
used in practice for a large variety of applications. Moreover, most
recent Machine Learning competitions such as the Pascal
VOC or Netflix challenge have been won by procedures combining
different types of
classifiers/predictors/estimators. It is therefore of central interest
to understand from a theoretical
point of view what kind of aggregation strategies should be used for
getting the best possible combination
of the available statistical procedures.

\subsection{Historical remarks and motivation}
In the statistical literature, to the best of our knowledge,
theoretical foundations of aggregation procedures
were first studied by Nemirovski (Nemirovski~\cite{Nemirovski00}, Juditsky and
Nemirovski~\cite{JuditskyNemirovski00})
and independently by a series of papers by Catoni (see~\cite{Catoni04}
for an account) and Yang~\cite{Yang00,Yang03,Yang04a}.
For the regression model, a significant progress was achieved by
Tsybakov~\cite{Tsybakov03} with introducing the notion of optimal
rates of aggregation and proposing aggregation-rate-optimal procedures
for the tasks of linear, convex and model selection
aggregation. This point was further developed in \cite
{Lounici07,RigolletTsybakov07,BuneaTsybakovWegkamp07},
especially in the context of high dimension with sparsity constraints
and in~\cite{Rigollet09} for Kullback--Leibler
aggregation. However, it should be noted that the procedures proposed
in~\cite{Tsybakov03} that provably achieve the lower bounds
in convex and linear aggregation require full knowledge of design
distribution. This limitation was overcome in the recent
work~\cite{Wang2011}.

From a practical point of view, an important limitation of the
previously cited results on aggregation is that
they are valid under the assumption that the aggregated procedures are
deterministic (or random, but independent
of the data used for aggregation). The generality of those
results---almost no restriction on the constituent
estimators---compensates to this practical limitation.

In the Gaussian sequence model, a breakthrough was reached by Leung and
Barron~\cite{LeungBarron06}. Building on very
elegant but not very well-known results by George~\cite{George86a}\setcounter{footnote}{1}\footnote{Corollary 2
in~\cite{George86a} coincides with Theorem 1 in
\cite{LeungBarron06}
in the case of exponential weights with temperature $\beta=2\sigma^2$;
cf. equation (\ref{eqdefweights}) below for a precise
definition of exponential weights.
Furthermore, to the best of our knowledge,~\cite{George86a} is the
first reference using the Stein lemma for evaluating the expected risk
of the exponentially
weighted aggregate.}, they established sharp oracle inequalities for
the exponentially weighted aggregate (EWA) for constituent estimators
obtained from
the data vector by orthogonally projecting it on some linear subspaces.
Dalalyan and
Tsybakov~\cite{DalalyanTsybakov07,DalalyanTsybakov08}
showed the result of~\cite{LeungBarron06} remains valid under more
general (non-Gaussian) noise distributions
and when the constituent estimators are independent of the data used
for the aggregation. A natural question arises
whether a similar result can be proved for a larger family of
constituent estimators containing projection estimators
and deterministic ones as specific examples. The main aim of the
present paper is to answer this question by considering
families of affine estimators.
%

Our interest in affine estimators is motivated by several reasons.
First, affine estimators encompass
many popular estimators such as smoothing splines, the Pinsker
estimator~\cite{Pinsker80,EfromovichPinsker96},
local polynomial estimators, nonlocal means \cite
{BuadesCollMorel05,SalmonLepennec09b}, etc. For instance,
it is known that if the underlying (unobserved) signal belongs to a
Sobolev ball, then the (linear) Pinsker
estimator is asymptotically minimax up to the optimal constant, while
the best projection estimator is only
rate-minimax. A second motivation is that---as proved by Juditsky and
Nemirovski~\cite{JuditskyNemirovski09}---the set of signals
that are well estimated by linear estimators is very rich. It contains,
for instance, sampled smooth functions,
sampled modulated smooth functions and sampled harmonic functions. One
can add to this set the family of piecewise
constant functions as well, as demonstrated in~\cite
{PolzehlSpokoiny00}, with natural application in magnetic
resonance imaging. It is worth noting that oracle inequalities for
penalized empirical risk minimizer were also
proved by Golubev~\cite{Golubev10}, and for model selection by Arlot
and Bach~\cite{ArlotBach09}, Baraud, Giraud and Huet \cite
{BaraudGiraudHuet10}.

In the present work, we establish sharp oracle inequalities in the
model of heteroscedastic regression, under various conditions
on the constituent estimators assumed to be affine functions of the data.
Our results provide theoretical guarantees of optimality, in terms of
expected loss, for the exponentially weighted aggregate.
They have the advantage of covering in a unified fashion the particular
cases of frozen estimators considered in~\cite{DalalyanTsybakov08}
and of projection estimators treated in~\cite{LeungBarron06}.

We focus on the theoretical guarantees expressed in terms of oracle
inequalities for
the expected squared loss. Interestingly, although several recent papers
\cite{ArlotBach09,BaraudGiraudHuet10,GoldenshlugerLepski08}
discuss the paradigm of
competing against the best linear procedure from a given family, none
of them provide oracle inequalities
with leading constant equal to one. Furthermore, most existing results
involve some constants depending
on different parameters of the setup. In contrast, the oracle
inequality that we prove
herein is with leading constant one and admits a simple formulation. It
is established for
(suitably symmetrized, if necessary) exponentially weighted aggregates~\cite{George86a,Catoni04,DalalyanTsybakov07}
with an arbitrary prior and a temperature parameter which is not too
small. The result is
nonasymptotic but leads to an asymptotically optimal residual term when
the sample size, as well
as the cardinality of the family of constituent estimators, tends to infinity.
In its general form, the residual term is similar to those obtained in the
PAC-Bayes setting~\cite{Mcallester98,Langford02,Seeger03} in that it
is proportional to the
Kullback--Leibler divergence between two probability distributions.

The problem of competing against the best procedure in a given family
was extensively studied in the
context of online learning and prediction with expert advice~\cite
{KivinenWarmuth99,Cesa-BianchiLugosi06}.
A connection between the results on online learning and statistical
oracle inequalities was established by Gerchinovitz~\cite{Gerchinovitz11}.

\subsection{Notation and examples of linear estimators}

Throughout this work, we focus on the heteroscedastic regression model
with Gaussian additive noise.
We assume we are given a vector
$\YY=(y_1,\ldots,y_n)^\top \in\R^n$ obeying the model
%
\begin{equation}
\label{eqmodel} y_i=f_i+\xi_i\qquad
\mbox{for } i=1,\ldots,n,
\end{equation}
where $\xxi=(\xi_1,\ldots, \xi_n)^\top $ is a centered Gaussian
random vector, $f_i=\mathbf{f}(x_i)$ where
$\mathbf{f}\dvtx\mathcal{X} \rightarrow\R$ is an unknown function and
$x_1,\ldots, x_n \in\mathcal{X}$ are
deterministic points. Here, no assumption is made on the set $\mathcal
{X}$. Our objective is to
recover the vector $\ff=(f_1,\ldots,f_n)^\top $, often referred to as
\textit{signal}, based on the data
$y_1,\ldots,y_n$. In our work, the noise covariance matrix $\Sigma=\E
[\xxi\xxi^\top ]$ is assumed to be
finite with a known upper bound on its spectral norm $|\!|\!|{\Sigma
}|\!|\!|$. We denote
by $ \langle\cdot| \cdot\rangle_n$ the empirical inner
product in $\R^n$:
$ \langle\bu| \bv\rangle_n=(1/n)\sum_{i=1}^n u_iv_i$.
We measure the performance of an estimator
$\hat{\mathbf{f}}$ by its expected empirical quadratic loss:
$r=\E[\|\ff-\hat{\mathbf{f}}\|_n^2 ]$ where $\|\ff-\hat{\mathbf
{f}}\|_n^2=\frac
{1}{n} \sum_{i=1}^{n} (f_i-\hat{f}_i)^2$.

We only focus on the task of aggregating \textit{affine estimators}
$\hat{\mathbf{f}}_{\lambda}$ indexed by some parameter $\lambda\in
\Lambda$.
These estimators can be written as affine transforms of the data $\YY
=(y_1,\ldots,y_n)^\top \in\R^n$.
Using the convention that all vectors are one-column matrices, we have
$\hat{\mathbf{f}}_{\lambda}=A_{\lambda}\YY+\mathbf{b}_{\lambda}$,
where the $n \times n$ real matrix $A_{\lambda}$ and the vector
$\mathbf{b}_{\lambda}\in
\R^n$ are deterministic. It means
the entries of $A_{\lambda}$ and $\mathbf{b}_{\lambda}$ may depend
on the points
$x_1,\ldots, x_n$ but not on the data $\YY$.
Let us describe now different families of linear and affine estimators
successfully used in the
statistical literature. Our results apply to all these families,
leading to a procedure that behaves nearly
as well as the best (unknown) one of the family.

\textit{Ordinary least squares}. Let $\{\mathcal S_\lambda\dvtx
\lambda
\in\Lambda\}$ be a set of linear subspaces
of $\R^n$.
A~well-known family of affine estimators, successfully used in the
context of model selection
\cite{BarronBirgeMassart99}, is the set of orthogonal
projections onto $\mathcal S_\lambda$. In the case of a family of
linear regression models with design
matrices $X_\lambda$, one has $A_\lambda=X_\lambda(X_\lambda^\top
X_\lambda)^+ X_\lambda^\top $, where
$(X_\lambda^\top X_\lambda)^+$ stands for the Moore--Penrose
pseudo-inverse\vspace*{1pt} of $X_\lambda^\top X_\lambda$.

\textit{Diagonal filters}. Other common estimators are
the so-called diagonal filters corresponding to diagonal matrices
$A=\diag(a_1,\ldots,a_n)$.
Examples include the following:

\begin{itemize}
\item Ordered projections: $a_k=\1_{(k \leq\lambda)}$
for some integer $\lambda$ [$\1_{(\cdot)}$ is
the indicator function]. Those weights are also called
truncated SVD (Singular Value Decomposition) or spectral cutoff.
In this case a natural parametrization is $\Lambda=\{1,\ldots, n\}$,
indexing the number of elements conserved.
\item Block projections: $a_k=\1_{(k\leq w_1 )} + \sum_{j=1}^{m-1}\lambda_j \1_{(w_j\leq k\leq w_{j+1})} $,
$k=1,\ldots,n$, where $\lambda_j\in\{0,1\}$. Here the natural
parametrization is $\Lambda= \{0,1\}^{m-1}$, indexing subsets of $\{
1,\ldots,m-1\}$.
\item Tikhonov--Philipps filter: $a_k=\frac
{1}{1+(k/w)^\alpha}$, where $w,\alpha>0$. In this case, $\Lambda=
(\R_+^{*})^2$, indexing continuously the smoothing parameters.
\item Pinsker filter: $a_k= (1-\frac{k^\alpha
}{w} )_{+}$, where $x_+=\max(x,0)$ and $(w,\alpha)=\lambda\in
\Lambda= (\R_+^{*})^2$.
\end{itemize}

\textit{Kernel ridge regression}.
Assume that we have a positive definite kernel $k\dvtx\mathcal
{X}\times
\mathcal{X} \rightarrow\R$ and we aim
at estimating the true function $f$ in the associated reproducing
kernel Hilbert space ($\mathcal H_k,\|\cdot\|_k$). The kernel ridge
estimator is obtained by minimizing
the criterion $\|\YY-\ff\|^2_n+\lambda\|\ff\|_k^2$ w.r.t. $f \in
\mathcal
H_k$ (see~\cite{Shawe-TaylorChristianini00}, page 118). Denoting by $K$ the $ n
\times n $ kernel-matrix with element
$K_{i,j}=k(x_i,x_j)$, the unique solution~$\hat{\mathbf{f}}$ is a
linear estimate
of the data, $\hat{\mathbf{f}}=A_{\lambda}\YY$, with
$A_{\lambda}=K(K+n\lambda{I_{n\times n}})^{-1}$, where ${I_{n\times
n}}$ is the $n\times n$
identity matrix.

\textit{Multiple Kernel learning}.
As described in~\cite{ArlotBach09}, it is possible to handle the case
of several kernels
$k_1,\ldots,k_M$, with associated positive definite matrices
$K_1,\ldots, K_M$.
For a parameter $\lambda=(\lambda_1,\ldots,\lambda_M) \in\Lambda=
\R_+^M$, one can define the estimators
$\hat{\mathbf{f}}_{\lambda}=A_\lambda\YY$ with
%
\begin{equation}
\label{eqmultiplekernel} A_\lambda= \Biggl(\sum
_{m=1}^M \lambda_m K_m
\Biggr) \Biggl(\sum_{m=1}^M
\lambda_m K_m + n {I_{n\times n}}
\Biggr)^{-1}.
\end{equation}
It is worth mentioning that the formulation in
equation \eqref{eqmultiplekernel} can be linked to the group Lasso
\cite{YuanLin06} and to
the multiple kernel learning introduced in \cite
{LanckrietCristianiniBartlettElGhaouiJordan03}---see
\cite{ArlotBach09} for more details.

\textit{Moving averages}. If we think of coordinates of $\ff$ as
some values assigned to the vertices of an
undirected graph, satisfying the property that two nodes are connected
if the corresponding values of $\ff$ are close, then it
is natural to estimate $f_i$ by averaging out the values $Y_j$ for
indices $j$ that are connected to $i$. The resulting estimator is
a linear one with a matrix $A=(a_{ij})_{i,j=1}^n$ such that $a_{ij}=\1_{V_i}(j)/n_i$, where $V_i$ is the set of neighbors of the node
$i$ in the graph and $n_i$ is the cardinality of $V_i$.

\subsection{Organization of the paper}
In Section~\ref{secmainresult} we introduce EWA and state a
PAC-Bayes type bound in expectation assessing optimality
properties of EWA in combining affine estimators. The strengths and
limitations of the results are discussed in Section~\ref{secdiscussion}.
The extension of these results to the case of grouped aggregation---in
relation with ill-posed inverse problems---is developed in Section~\ref
{secill-posed}. As a consequence, we provide in
Section~\ref{secexamples} sharp oracle inequalities in various
setups: ranging from finite to continuous families of
constituent estimators and including sparse scenarii. In Section \ref
{secmmx} we apply our main results to prove
that combining Pinsker's type filters with EWA leads to asymptotically
sharp adaptive procedures over Sobolev ellipsoids.
Section~\ref{secexperiments} is devoted to numerical comparison of
EWA with other classical filters
(soft thresholding, blockwise shrinking, etc.) and illustrates
the potential benefits of aggregating. The conclusion is given in
Section~\ref{secconclusion}, while
the proofs of some technical results (Propositions~\ref{prop2}--\ref{propminimax2}) are provided in
the supplementary material~\cite{DalalyanSalmonsupp}.

\section{Aggregation of estimators: Main results}\label{secmainresult}
In this section we describe the statistical framework for aggregating estimators
and we introduce the exponentially\vadjust{\goodbreak} weighted aggregate. The task of aggregation
consists in estimating $\ff$ by a suitable combination of the elements
of a family
of \textit{constituent estimators} $\mathcal F_\Lambda=(\hat{\mathbf
{f}}_{\lambda}
)_{\lambda\in\Lambda}
\in\R^n$. The target objective of the aggregation is to build an aggregate
$\hat{\mathbf{f}}_{\mathrm{aggr}}$ that mimics the performance of
the best
constituent estimator,
called \textit{oracle} (because of its dependence on the unknown
function $\mathbf{f}$).
In what follows, we assume that $\Lambda$ is a measurable subset of
$\R^M$, for
some $M\in\mathbb N$.

The theoretical tool commonly used for evaluating the quality of an aggregation
procedure is the oracle inequality (OI), generally written
%
\begin{equation}
\label{ineqoracleex} \E \bigl[\|\hat{\mathbf{f}}_{\mathrm{aggr}} -\ff
\|^2_ n \bigr] \leq C_n \inf_{\lambda\in\Lambda} \E
\bigl[\|\hat{\mathbf {f}}_{\lambda}-\ff \|^2_ n
\bigr] +R_n,
\end{equation}
with \textit{residual} term $R_n$ tending to zero as $n\to\infty$,
and \textit{leading constant}
$C_n$ being bounded. The OIs with leading constant one are of central
theoretical
interest since they allow to bound the excess risk and to assess the
aggregation-rate-optimality. They are often referred to as sharp OI.

\subsection{Exponentially weighted aggregate (EWA)}

Let $r_\lambda=\E[\|\hat{\mathbf{f}}_{\lambda}-\ff\|_n^2 ]$
denote the risk
of the estimator $\hat{\mathbf{f}}_{\lambda}$, for any
$\lambda\in\Lambda$, and let $\hat{r}_{\lambda}$ be an estimator of
$r_\lambda$. The precise form of
$\hat{r}_{\lambda}$ strongly depends on the nature of the constituent
estimators. For any probability
distribution $\pi$ over $\Lambda$ and for any $\beta>0$, we define
the probability measure of
exponential weights, $\hat\pi$, by
%
\begin{equation}
\label{eqdefweights} \hat\pi(d\lambda)=\theta(\lambda)\pi (d\lambda) \qquad\mbox
{with } \theta(\lambda)=\frac{\exp(-n\hat{r}_{\lambda}/\beta)
}{\int_{\Lambda}\exp
(-n\hat{r}_{\omega}/\beta) \pi(d\omega) }.
\end{equation}
The corresponding exponentially weighted aggregate, henceforth denoted
by ${\hat{\mathbf{f}}}_\mathrm{EWA}$, is the expectation
of $\hat{\mathbf{f}}_{\lambda}$ w.r.t.  the probability
measure $\hat{\pi}$:
%
\begin{equation}
\label{eqdefestimator} {\hat{\mathbf{f}}}_\mathrm{EWA}= \int
_{\Lambda}\hat{\mathbf{f}}_{\lambda}~\hat{\pi}(d\lambda) .
\end{equation}
We will frequently use the terminology of Bayesian statistics: the
measure $\pi$ is called
\textit{prior}, the measure $\hat\pi$ is called
\textit{posterior} and the aggregate ${\hat{\mathbf{f}}}_\mathrm
{EWA}$ is then the \textit
{posterior mean}. The parameter $\beta$
will be referred to as the \textit{temperature parameter}. In the
framework of aggregating statistical
procedures, the use of such an aggregate can be traced back to
George~\cite{George86a}.

The interpretation of the weights $\theta(\lambda)$ is simple: they
up-weight estimators all the more that
their performance, measured in terms of the risk estimate $\hat
{r}_{\lambda}$,
is good. The temperature parameter reflects the confidence we have in
this criterion: if the temperature is
small ($\beta\approx0$), the distribution concentrates
on the estimators achieving the smallest value for $\hat{r}_{\lambda}$,
assigning almost zero weights to the other
estimators. On the other hand, if $\beta\rightarrow+\infty$, then
the probability distribution over $\Lambda$
is simply the prior $\pi$, and the data do not influence our
confidence in the estimators.

\subsection{Main results}

In this paper we only focus on \textit{affine estimators}
%
\begin{equation}
\label{eqaffineestimatorsgal} \hat{\mathbf{f}}_{\lambda}=A_{\lambda
}
\YY+\mathbf{b}_{\lambda},
\end{equation}
where the $n \times n$ real matrix $A_{\lambda}$ and the vector
$\mathbf{b}_{\lambda}\in
\R^n$ are deterministic.
Furthermore, we will assume that an unbiased estimator $\widehat\Sigma
$ of the
noise covariance matrix~$\Sigma$ is available.
It is well known (\lcf \hyperref[secproofs]{Appendix} for details) that
the risk of the
estimator \eqref{eqaffineestimatorsgal} is given by
%
\begin{equation}
\label{eqaffineestimatorsrisk} r_\lambda=\E \bigl[\|\hat{\mathbf
{f}}_{\lambda}-\ff \|_n^2 \bigr]=
\bigl\|(A_{\lambda}-I_{n\times n}) \ff+\mathbf{b}_{\lambda
}
\bigr\|_n^2+ \frac{\Tr(A_{\lambda}\Sigma A_{\lambda}^\top )}{n}
\end{equation}
and that $\hat r_\lambda^{\mathrm{unb}}$, defined by
%
\begin{equation}
\label{equnbiasedrisk} \hat r_\lambda^{\mathrm{unb}}= \|\YY- \hat{
\mathbf{f}}_{\lambda} \|^2_n+\frac{2}{n} \Tr(
\widehat \Sigma A_{\lambda})-\frac{1}{n}\Tr[\widehat\Sigma],
\end{equation}
is an unbiased estimator of $r_\lambda$. Along with $\hat r_\lambda^{\mathrm{unb}}$, we
will use another estimator of the risk
that we call the adjusted risk estimate and define by
%
\begin{equation}
\label{equnbiasedrisk1} \hat r_\lambda^{\mathrm{adj}}=\underbrace { \|
\YY- \hat{\mathbf{f}}_{\lambda} \|^2_n+
\frac{2}{n} \Tr(\widehat\Sigma A_{\lambda})-\frac{1}{n}\Tr[
\widehat\Sigma ]}_{\hat r_\lambda^{\mathrm{unb}}} + \frac{1}{n} \YY^\top
\bigl(A_\lambda-A_\lambda^2 \bigr)\YY.
\end{equation}
One can notice that the adjusted risk estimate $\hat r_\lambda^{\mathrm{adj}}$ coincides
with the unbiased risk estimate
$\hat r_\lambda^{\mathrm{unb}}$ if and only if the matrix $A_{\lambda
}$ is an orthogonal projector.

To state our main results, we denote by $\mathcal{P}_\Lambda$ the set
of all
probability measures on $\Lambda$ and by
$\KL(p,p')$ the Kullback--Leibler divergence
between two probability measures $p,p'\in\mathcal{P}_\Lambda$:
\[
\mathcal{K} \bigl(p,p' \bigr)=\cases{ \displaystyle\int_{\Lambda}
\log \biggl( \frac{dp}{dp'} (\lambda) \biggr)p(d\lambda), & \quad $\mbox{if } p
\mbox{ is absolutely continuous w.r.t. } p',$\vspace*{2pt}
\cr
+\infty,& \quad $\mbox{otherwise.}$}
\]
We write $S_1 \preceq S_2$ (resp., $S_1 \succeq S_2$) for two symmetric
matrices $S_1$ and $S_2$, when $S_2 -S_1$ (resp., $S_1-S_2$)
is semi-definite positive.

%
\begin{theorem}\label{mainthm}
Let all the matrices $A_{\lambda}$ be symmetric and $\widehat\Sigma$ be
unbiased and independent of $\YY$.
\begin{longlist}[(ii)]
\item[(i)] Assume that for all $\lambda,\lambda'\in\Lambda$,
it holds that $A_{\lambda}A_{\lambda'}=A_{\lambda'}A_{\lambda}$,
$A_{\lambda}\Sigma
+\Sigma A_{\lambda}\succeq0$
and $\mathbf{b}_{\lambda}=0$. If $\beta\ge8|\!|\!|{\Sigma}|\!|\!
|$, then the aggregate
${\hat{\mathbf{f}}}_\mathrm{EWA}$
defined by equations (\ref{eqdefweights}), (\ref
{eqdefestimator}) and the unbiased risk estimate $\hat{r}_{\lambda}
=\hat r_\lambda^{\mathrm{unb}}$ (\ref{equnbiasedrisk}) satisfies
%
\begin{equation}
\label{ineqKL} \E \bigl[\|{\hat{\mathbf{f}}}_\mathrm{EWA}-\ff
\|_n^2 \bigr] \leq \inf_{p \in\mathcal{P}_\Lambda} \biggl\{ \int
_{\Lambda}\E \bigl[\|\hat{\mathbf{f}}_{\lambda}-\ff
\|_n^2 \bigr] p(d\lambda)+ \frac{\beta}n \KL(p,\pi)
\biggr\}.
\end{equation}
\item[(ii)] Assume that, for all $\lambda\in\Lambda$,
$A_{\lambda}
\preceq I_{n\times n}$ and $A_{\lambda}\mathbf{b}_{\lambda}=0$.
If $\beta\ge4|\!|\!|{\Sigma}|\!|\!|$, then the aggregate ${\hat
{\mathbf{f}}}_\mathrm{EWA}$
defined by equations (\ref{eqdefweights}), (\ref{eqdefestimator})
and the adjusted risk estimate $\hat{r}_{\lambda}=\hat r_\lambda^{\mathrm{adj}}$ (\ref
{equnbiasedrisk1}) satisfies
\begin{eqnarray*}
\E \bigl[\|{\hat{\mathbf{f}}}_\mathrm{EWA}-\ff
\|_n^2 \bigr] &\leq& \inf_{p \in\mathcal{P}_\Lambda} \biggl\{ \int
_{\Lambda}\E \bigl[\|\hat{\mathbf{f}}_{\lambda}-\ff
\|_n^2 \bigr] p(d\lambda)+ \frac{\beta}n \KL(p,\pi)
\nonumber
\\
&&\hspace*{26pt}{} +\frac1n\int_{\Lambda} \bigl(\ff^\top
\bigl(A_{\lambda
}-A_{\lambda}^2 \bigr) \ff+\Tr \bigl[\Sigma
\bigl(A_{\lambda}-A_{\lambda}^2 \bigr) \bigr] \bigr)p(d
\lambda) \biggr\}.
\end{eqnarray*}
\end{longlist}
\end{theorem}

The simplest setting in which all the conditions of part (i) of
Theorem~\ref{mainthm} are fulfilled is when the matrices $A_{\lambda}$
and $\Sigma$ are all diagonal, or diagonalizable in a common base.
This result, as we will see in Section~\ref{secmmx},
leads to a new estimator which is adaptive, in the exact minimax sense,
over the collection of all Sobolev ellipsoids. It also
suggests a new method for efficiently combining varying-block-shrinkage
estimators, as described in Section~\ref{sec2stein}.

However, part (i) of Theorem~\ref{mainthm} leaves open the issue of
aggregating affine estimators defined via noncommuting matrices.
In particular, it does not allow us to evaluate the MSE of EWA when
each $A_{\lambda}$ is a convex or linear combination
of a fixed family of projection matrices on nonorthogonal linear
subspaces. These kinds of situations may be handled via
the result of part (ii) of Theorem~\ref{mainthm}. One can observe
that in the particular case of a finite collection of
projection estimators (i.e., $A_{\lambda}=A_{\lambda}^2$ and $\mathbf
{b}_{\lambda}=0$
for every $\lambda$), the result of part (ii) offers
an extension of~\cite{LeungBarron06}, Corollary 6, to the case of
general noise covariances (\cite{LeungBarron06} deals only
with i.i.d. noise).

An important situation covered by part (ii) of Theorem~\ref{mainthm},
but not by part (i), concerns the case when signals of
interest $\ff$ are smooth or sparse in a basis $\mathcal{B}_{\mathrm
{sig}}$ which is
different from the basis $\mathcal{B}_{\mathrm{noise}}$ orthogonalizing
the covariance matrix $\Sigma$. In such a context, one may be
interested in considering matrices $A_{\lambda}$ that are diagonalizable
in the basis $\mathcal{B}_{\mathrm{sig}}$ which, in general, do not
commute with $\Sigma$.

%
\begin{remark}
While the results in~\cite{LeungBarron06} yield a sharp oracle
inequality in the case of projection matrices $A_{\lambda}$, they are of
no help in the case when the matrices $A_{\lambda}$ are nearly idempotent
and not exactly. Assertion (ii) of Theorem~\ref{mainthm}
fills this gap by showing that if $\max_\lambda\Tr[A_{\lambda
}-A_{\lambda}^2]\le\delta$, then $\E[\|{\hat{\mathbf{f}}}_\mathrm
{EWA}-\ff\|_n^2 ]$ is bounded
by
\[
\inf_{p \in\mathcal{P}_\Lambda} \biggl\{ \int_{\Lambda}\E \bigl[\|\hat{
\mathbf{f}}_{\lambda}-\ff\|_n^2 \bigr] p(d\lambda)+
\frac{\beta}n \KL(p,\pi) \biggr\}+\delta \bigl(\|\ff \|_n^2+n^{-1}|
\!|\!|{\Sigma}|\!|\!| \bigr).
\]
\end{remark}

%
\begin{remark}
We have focused only on Gaussian errors to emphasize that it is
possible to efficiently aggregate almost any family
of \textit{affine estimators}. We believe that by a suitable
adaptation of the approach developed in~\cite{DalalyanTsybakov08},
claims of Theorem~\ref{mainthm} can be generalized---at least when
$\xi_i$ are independent with known variances---to some
other common noise distributions.
\end{remark}

The results presented so far concern the situation when the matrices
$A_{\lambda}$ are symmetric. However, using
the last part of Theorem~\ref{mainthm}, it is possible to propose an
estimator of $\ff$ that is almost as
accurate as the best affine estimator $A_{\lambda}\YY+\mathbf
{b}_{\lambda}$ even if the
matrices $A_{\lambda}$ are not symmetric.
Interestingly, the estimator enjoying this property is not obtained by
aggregating the original\vspace*{1pt} estimators
$\hat{\mathbf{f}}_{\lambda}=A_{\lambda}\YY+\mathbf{b}_{\lambda}$
but the ``symmetrized'' estimators
$\tilde{\ff}_{\lambda}=\tilde{A}_{\lambda}\YY+\mathbf
{b}_{\lambda}$, where
$\tilde{A}_{\lambda}=A_{\lambda}+A_{\lambda}^\top -A_{\lambda
}^\top A_{\lambda}$. Besides symmetry, an
advantage of the matrices $\tilde{A}_{\lambda}$, as
compared to the $A_{\lambda}$'s, is that they automatically satisfy
the contraction
condition $\tilde{A}_{\lambda}\preceq I_{n\times n}$ required by part
(ii) of
Theorem~\ref{mainthm}.
We will refer to this method as Symmetrized Exponentially Weighted
Aggregates (or SEWA)
\cite{DalalyanSalmon11b}.

%
\begin{theorem}\label{mainthmcor}
Assume that the matrices $A_{\lambda}$ and the vectors $\mathbf
{b}_{\lambda}$ satisfy
$A_{\lambda}\mathbf{b}_{\lambda}=A_{\lambda}^\top \mathbf
{b}_{\lambda}=0$ for every $\lambda\in\Lambda$.
Assume in addition that $\widehat\Sigma$ is an unbiased estimator of
$\Sigma$ and is independent of $\YY$. Let
$\tilde{\mathbf{f}}_\mathrm{SEWA}$ denote the exponentially weighted
aggregate of the
(symmetrized) estimators $\tilde{\ff}_{\lambda}=(A_{\lambda
}+A_{\lambda}^\top -A_{\lambda}^\top A_{\lambda})\YY
+\mathbf{b}_{\lambda}$
with the weights~(\ref{eqdefweights}) defined via the risk estimate
$\hat r_\lambda^{\mathrm{unb}}$. Then, under the conditions
$\beta\ge4|\!|\!|{\Sigma}|\!|\!|$ and
%
{\renewcommand{\theequation}{\Alph{equation}}
\setcounter{equation}{2}
\begin{equation}\label{eqC}
\pi \biggl\{\lambda\in\Lambda\dvtx\Tr (\widehat\Sigma
A_{\lambda})\le\Tr\bigl(\widehat\Sigma A_{\lambda}^\top
A_{\lambda}\bigr) \biggr\}=1 \qquad\mbox{a.s.}
\end{equation}}
\hspace*{-2pt}it holds that
\renewcommand{\theequation}{\arabic{section}.\arabic{equation}}
\setcounter{equation}{8}
\begin{equation}
\label{ineqKL4} \E \bigl[\|\tilde{\mathbf{f}}_\mathrm{SEWA}-\ff
\|_n^2 \bigr]  \leq\inf_{p \in\mathcal{P}_\Lambda} \biggl\{ \int
_{\Lambda}\E \bigl[\|\hat{\mathbf{f}}_{\lambda}-\ff
\|_n^2 \bigr]p(d\lambda)+ \frac{\beta}n \KL(p,\pi)
\biggr\}.
\end{equation}
\end{theorem}

To understand the scope of condition {(\ref{eqC})}, let us present
several cases of widely used linear estimators for which this
condition is satisfied:
\begin{itemize}
\item The simplest class of matrices $A_{\lambda}$ for which condition
(\ref{eqC}) holds true are orthogonal projections.
Indeed, if $A_{\lambda}$ is a projection matrix, it satisfies
$A_{\lambda}^\top A_{\lambda}
=A_{\lambda}$ and, therefore, $\Tr(\widehat\Sigma A_{\lambda})=
\Tr(\widehat\Sigma A_{\lambda}^\top A_{\lambda})$.
\item When the matrix $\widehat\Sigma$ is diagonal, then a sufficient
condition for (\ref{eqC}) is $a_{ii}\le\sum_{j=1}^n a_{ji}^2$.
Consequently, (\ref{eqC}) holds true for matrices having only zeros
on the main diagonal. For instance, the $k$NN filter
in which the weight of the observation $Y_i$ is replaced by zero,
that is, $a_{ij}=\mathbf{1}_{j\in\{j_{i,1},\ldots,j_{i,k}\}}/k$
satisfies this condition.
\item Under a little bit more stringent assumption of homoscedasticity,
that is, when $\widehat\Sigma=\widehat\sigma^2I_{n\times n}$,
if the matrices $A_{\lambda}$ are such that all the nonzero elements
of each
row are equal and sum up to one (or a
quantity larger than one), then $\Tr(A_{\lambda})=\Tr(A_{\lambda
}^\top A_{\lambda})$ and
(\ref{eqC}) is fulfilled. A notable
example of linear estimators that satisfy this condition are
Nadaraya--Watson estimators with rectangular kernel
and nearest neighbor filters.
\end{itemize}

\section{Discussion}\label{secdiscussion}

Before elaborating on the main results stated in the previous section,
by extending them to inverse problems and by
deriving adaptive procedures, let us discuss some aspects of the
presented OIs.

\subsection{\texorpdfstring{Assumptions on $\Sigma$}{Assumptions on Sigma}}
In some rare situations, the matrix $\Sigma$ is known and it is natural to use
$\widehat\Sigma=\Sigma$ as an unbiased estimator. Besides this not
very realistic situation, there are at least two
contexts in which it is reasonable to assume that an unbiased estimator
of $\Sigma$, independent of $\YY$, is available.

The first case corresponds to problems in which a signal can be
recorded several times by the same device, or once
but by several identical devices. For instance, this is the case when
an object is photographed many times by
the same digital camera during a short time period. Let $\bZ_1,\ldots
,\bZ_N$ be the available signals, which
can be considered as i.i.d. copies of an $n$-dimensional Gaussian
vector with mean $\ff$ and covariance
matrix $\Sigma_Z$. Then, defining $\YY=(\bZ_1+\cdots+\bZ_N)/N$ and
$\widehat\Sigma_Z=(N-1)^{-1}(\bZ_1\bZ_1^\top +\cdots+\bZ_N\bZ_N^\top
-N\YY\YY^\top )$,
we find ourselves within the framework covered by previous theorems.
Indeed, $\YY\sim\mathcal N_n(\ff,\Sigma_Y)$
with $\Sigma_Y=\Sigma_Z/N$ and $\widehat\Sigma_Y=\widehat\Sigma_Z/N$ is an unbiased estimate of $\Sigma_Y$, independent
of $\YY$. Note that our theory applies in this setting for every
integer $N\ge2$.

The second case is when the dominating part of the noise comes from the
device which is used for recording the signal.
In this case, the practitioner can use the device in order to record a
known signal, $\boldg$. In digital image
processing, $\boldg$ can be a black picture. This will provide a noisy
signal $\bZ$ drawn from Gaussian distribution
$\mathcal N_n(\boldsymbol{g},\Sigma)$, independent of $\YY$ which is
the signal of interest. Setting
$\widehat\Sigma=(\bZ-\boldg)(\bZ-\boldg)^\top $, one ends up with
an unbiased estimator of $\Sigma$, which is
independent of $\YY$.

\subsection{OI in expectation versus OI with high probability}
All the results stated in this work provide sharp nonasymptotic bounds
on the expected risk of EWA. It would
be insightful to complement this study by risk bounds that hold true
with high probability. However, it was recently
proved in~\cite{DaiRigZha12} that EWA is deviation suboptimal: there
exist a family of constituent estimators and a
constant $C>0$ such that the difference between the risk of EWA and
that of the best constituent estimator is larger
than $C/\sqrt{n}$ with probability at least $0.06$. Nevertheless,
several empirical studies (see, e.g.,~\cite{DaiZhang11})
demonstrated that EWA has often a smaller risk than some of its
competitors, such as the empirical star procedure~\cite{Audibert07},
which are provably optimal in the sense of OIs with high probability.
Furthermore, numerical experiments carried out in
Section~\ref{secexperiments} show that the standard-deviation of the
risk of EWA is of the order of $1/n$. This
suggests that under some conditions on the constituent estimators it
might be possible to establish OIs for EWA that
are similar to (\ref{ineqKL}) but hold true with high probability. A
step in proving this kind of result
was done in~\cite{LecueMendelson10}, Theorem~C, for the model of
regression with random design.

\subsection{Relation to previous work and limits of our results}
The OI of the previous section requires various conditions on the
constituent estimators
$\hat{\mathbf{f}}_{\lambda}=A_{\lambda}\YY+\mathbf{b}_{\lambda
}$. One may wonder how general these
conditions are and is it possible to extend
these OIs to more general $\hat{\mathbf{f}}_{\lambda}$'s. Although
this work does not
answer this question, we
can sketch some elements of response.

First of all, we stress that the conditions of the present paper relax
significantly those of previous
results existing in statistical literature. For instance, Kneip \cite
{Kneip94} considered only linear estimators,
that is, $\mathbf{b}_{\lambda}\equiv0$ and, more importantly, only
ordered sets of \textit{commuting} matrices $A_{\lambda}$.
The ordering assumption is dropped in~Leung and Barron \cite
{LeungBarron06}, in the case of projection matrices. Note that
neither of these assumptions is satisfied for the families of Pinsker
and Tikhonov--Philipps estimators. The
present work strengthens existing results in considering more general,
affine estimators extending both
projection matrices and ordered commuting matrices.

Despite the advances achieved in this work, there are still interesting
cases that are not covered by
our theory. We now introduce a family of estimators commonly used in
image processing that do not
satisfy our assumptions. In recent years, nonlocal means (NLM) became
quite popular in image
processing~\cite{BuadesCollMorel05}. This method of signal
denoising, shown to be tied in with EWA~\cite{SalmonLepennec09b}, removes noise by exploiting signals
self-similarities. We briefly define the
NLM procedure in the case of one-dimensional signals.

Assume that a vector $\YY=(y_1,\ldots,y_n)^\top $ given by (\ref
{eqmodel}) is observed with $f_i=\mathbf{f}(i/n)$,
$i=1,\ldots,n$, for some function $\mathbf{f}\dvtx[0,1]\to\R$. For
a fixed
``patch-size'' $k \in\{1,\ldots,n\}$, let us define $\ff_{[i]}=(f_i,f_{i+1},\ldots,f_{i+k-1})^\top $ and
$\YY_{ [i]}=(y_i,y_{i+1},\ldots,y_{i+k-1})^\top $ for every
$i=1,\ldots,n-k+1$. The vectors
$\ff_{[i]}$ and $\YY_{ [i]}$ are, respectively, called \textit{true
patch} and \textit{noisy patch}.
The NLM consists in regarding the noisy patches $\YY_{ [i]}$ as
constituent estimators for estimating
the true patch $\ff_{[i_0]}$ by applying EWA. One easily checks that
the constituent estimators
$\YY_{ [i]}$ are affine in $\YY_{ [i_0]}$, that is, $\YY_{
[i]}=A_i\YY_{ [i_0]}+\bb_i$ with $A_i$
and $\bb_i$ independent of $\YY_{ [i_0]}$. Indeed, if the distance
between $i$ and $i_0$ is larger
than $k$, then $\YY_{ [i]}$ is independent of $\YY_{ [i_0]}$ and,
therefore, $A_i=0$ and
$\bb_i=\YY_{ [i]}$. If $|i-i_0|<k$, then the matrix $A_i$ is a
suitably chosen shift matrix and $\bb_i$
is the projection of $\YY_{ [i]}$ onto the orthogonal complement of
the image of $A_i$. Unfortunately,
these matrices $\{A_i\}$ and vectors $\{b_i\}$ do not fit our
framework, that is, the assumption
$A_i b_i=A_i^\top b_i=0$ is not satisfied.

Finally, our proof technique is specific to affine estimators. Its
extension to estimators
defined as a more complex function of the data will certainly require
additional tools and is a challenging
problem for future research. Yet, it seems unlikely to get sharp OIs
with optimal remainder term for a fairly
general family of constituent estimators (without data-splitting),
since this generality inherently
increases the risk of overfitting.

\section{Ill-posed inverse problems and group-weighting}\label{secill-posed}

As explained in~\cite{CavalierGolubevPicardTsybakov02,Cavalier08},
the model of heteroscedastic regression
is well suited for describing inverse problems. In fact, let $T$ be a
known linear operator on some Hilbert space $\mathcal H$,
with inner product $ \langle\cdot| \cdot\rangle_{\mathcal H}$. For
some $h
\in\mathcal H$, let $Y$ be the random process
indexed by $g\in\mathcal H$ such that
%
\begin{equation}
\label{eqinversepb} Y=T h+\varepsilon\xi\quad\Longleftrightarrow \quad\bigl( Y(g)= \langle
T h | g \rangle_{\mathcal
H}+\varepsilon\xi(g), \forall g \in\mathcal H \bigr),
\end{equation}
where $\varepsilon>0$ is the noise magnitude and $\xi$ is a white
Gaussian noise on $\mathcal H$,
that is, for any $g_1,\ldots,g_k\in\mathcal H$
the vector $ (Y(g_1),\ldots,Y(g_k) )$ is Gaussian with zero
mean and covariance matrix $\{\langle
g_i | g_j\rangle_{\mathcal H}\}$. The problem is then the following:
estimate the element $h$ assuming the value of $Y$ can be measured for
any given $g$.
It is customary to use as $g$ the eigenvectors of the adjoint $T^*$ of $T$.
Under the condition that the operator $T^* T$ is compact, the SVD yields
$T \phi_k=b_k \psi_k$ and $T^* \psi_k=b_k \phi_k$, for $k\in
\mathbb N$,
where $b_k$ are the singular values, $\{\psi_k\}$ is an orthonormal
basis in $\Range(T) \subset\mathcal H$
and $\{\phi_k\}$ is the corresponding orthonormal basis in $\mathcal
H$. In view of~\eqref{eqinversepb}, it holds that
%
\begin{equation}
\label{eqinversepb2} Y(\psi_k)= \langle h | \phi_k
\rangle_{\mathcal H} b_k +\varepsilon\xi(\psi_k),\qquad k\in
\mathbb N.
\end{equation}
Since in practice only a finite number of measurements can be computed,
it is natural to assume that the
values $Y(\psi_k)$ are available
only for $k$ smaller than some integer $n$. Under the assumption that
$b_k \neq0$,\vspace*{1pt} the last equation is
equivalent to \eqref{eqmodel} with
$f_i= \langle h | \phi_i \rangle_{\mathcal H}$ and
$\Sigma=\diag(\sigma_i^2,i=1,2,\ldots)$ for $\sigma_i=\varepsilon
b_i^{-1}$.
Examples of inverse problems to which this statistical model has been
successfully applied are derivative estimation, deconvolution with
known kernel,
computerized tomography---see~\cite{Cavalier08} and the references
therein for more applications.

For very mildly ill-posed inverse problems, that is, when the
singular values $b_k$ of $T$
tend to zero not faster than any negative power of $k$, the approach
presented in Section~\ref{secmainresult}
will lead to satisfactory results. Indeed, by choosing $\beta
=8|\!|\!|{\Sigma}|\!|\!|$ or $\beta=4|\!|\!|{\Sigma}|\!|\!|$,
the remainder term in (\ref{ineqKL}) and (\ref{ineqKL4})
becomes---up to a logarithmic factor---proportional
to $\max_{1\le k\le n} b_k^{-2}/{n}$, which is the optimal rate in the
case of very mild ill-posedness.

However, even for mildly ill-posed inverse problems, the approach
developed in the previous section becomes
obsolete since the remainder blows up when $n$ increases to infinity.
Furthermore, this is not an artifact
of our theoretical results, but rather a drawback of the aggregation
strategy adopted in the previous section.
Indeed, the posterior probability measure $\hat\pi$ defined by (\ref
{eqdefweights}) can be seen as the
solution of the entropy-penalized empirical risk minimization problem:
%
\begin{equation}
\hat\pi_n=\operatorname{arg}\inf_{p} \biggl\{ \int
_\Lambda\hat{r}_{\lambda}p(d\lambda)+ \frac{\beta}{n}
\KL(p,\pi) \biggr\},
\end{equation}
where the $\inf$ is taken over the set of all probability
distributions. It means the same
regularization parameter $\beta$ is employed for estimating both the
coefficients $f_i= \langle
h | \phi_i \rangle_{\mathcal H}$ corrupted by noise of small
magnitude and those corrupted by
large noise. Since we place ourselves in the setting of known operator
$T$ and, therefore, known noise
levels, such a uniform treatment of all coefficients is unreasonable.
It is more natural to
upweight the regularization term in the case of large noise
downweighting the data fidelity term and,
conversely, to downweight the regularization in the case of small
noise. This motivates our interest in
the grouped EWA (or GEWA).

Let us consider a partition $B_1,\ldots,B_J$ of the set $\{1,\ldots
,n\}$: $B_j=\{T_j+1,\ldots,T_{j+1}\}$, for
some integers $0=T_1<T_2<\cdots<T_{J+1}=n$. To each element $B_j$ of this
partition, we associate the data sub-vector $\YY^{j}=(Y_i\dvtx i\in B_j)$
and the sub-vector of true function
$\ff^{j}=(f_i\dvtx i\in B_j)$. As in previous sections, we are
concerned by
the aggregation of affine estimators
$\hat{\mathbf{f}}_{\lambda}=A_{\lambda}\YY+\mathbf{b}_{\lambda
}$, but here we will assume the matrices
$A_{\lambda}$ are block-diagonal:
\[
A_{\lambda}=\left[\matrix{ A_{\lambda}^1 & 0 & \ldots& 0
\vspace*{2pt}
\cr
0 & A_{\lambda}^2 & \ldots& 0\vspace*{2pt}
\cr
\vdots& \vdots& \ddots& \vdots \vspace*{2pt}
\cr
0 & 0 & \ldots&
A_{\lambda}^J }\right] \qquad\mbox{with } A_{\lambda}^j\in
\R^{(T_{j+1}-T_j)\times
(T_{j+1}-T_j)}.
\]
Similarly, we define $\hat{\mathbf{f}}_{\lambda}^j$ and $\mathbf
{b}_{\lambda}^j$ as the sub-vectors of
$\hat{\mathbf{f}}_{\lambda}$ and $\mathbf{b}_{\lambda}$, respectively,
corresponding to the indices belonging to $B_j$. We will also assume
that the noise covariance matrix $\Sigma$ and its unbiased estimate
$\widehat\Sigma$ are block-diagonal with $(T_{j+1}-T_j)\times
(T_{j+1}-T_j)$ blocks $\Sigma^j$ and $\widehat\Sigma^j$,
respectively. This
notation implies, in particular, that $\hat{\mathbf{f}}_{\lambda
}^j=A_{\lambda}^j\YY^j+\mathbf{b}_{\lambda}^j$
for every $j=1,\ldots,J$. Moreover, the unbiased risk estimate
$\hat r_\lambda^{\mathrm{unb}}$ of $\hat{\mathbf{f}}_{\lambda}$
can be decomposed into the sum of unbiased
risk estimates $\hat r_\lambda^{j,\mathrm{unb}}$ of $\hat{\mathbf
{f}}_{\lambda}^j$, namely, $\hat r_\lambda^{\mathrm{unb}}
=\sum_{j=1}^J \hat r_\lambda^{j,\mathrm{unb}}$, where
\[
\hat r_\lambda^{j,\mathrm{unb}}=\bigl\|\YY^j-\hat{
\mathbf{f}}_{\lambda
}^j\bigr\|+\frac{2}{n} \Tr \bigl(\widehat
\Sigma^j A_{\lambda}^j \bigr)-\frac{1}{n}\Tr
\bigl[\widehat\Sigma^j \bigr], \qquad j=1,\ldots,J.
\]
To state the analogues of Theorems~\ref{mainthm} and~\ref
{mainthmcor}, we introduce the following settings.
\begin{longlist}
\item[\textit{Setting} 1:] For all
$\lambda,\lambda'\in\Lambda$ and $j\in\{1,\ldots,J\}$,
$A_{\lambda}^j$
are symmetric and satisfy $A_{\lambda}^j A_{\lambda'}^j=A_{\lambda
'}^jA_{\lambda}^j$,
$A_{\lambda}^j \Sigma^j+\Sigma^jA_{\lambda}^j\succeq0$ and
$\mathbf{b}_{\lambda}^j=0$. For a
temperature vector
$\bbeta=(\beta_1,\ldots,\beta_J)^\top $ and a\vadjust{\goodbreak} prior $\pi$, we
define GEWA
as ${\hat{\mathbf{f}}}_\mathrm{GEWA}^j=\int_\Lambda\hat{\mathbf
{f}}_{\lambda}^j \hat\pi^j(d\lambda)$,
where $\hat\pi^j(d\lambda)=\theta^j(\lambda)\pi(d\lambda)$ with
%
\begin{equation}
\label{eqdefgweights} \theta^j(\lambda)=\frac{\exp(-n\hat
r_\lambda^{j,\mathrm{unb}}/\beta_j) }{\int_{\Lambda}
\exp(-n\hat{r}_{\omega}^{j,\mathrm{unb}}/\beta_j) \pi(d\omega) }.
\end{equation}
\item[\textit{Setting} 2:] For every $j=1,\ldots,J$ and for every $\lambda$
belonging to a set of $\pi$-measure one, the matrices
$A_{\lambda}$ satisfy a.s. the inequality $\Tr(\widehat\Sigma^jA_{\lambda}^j)\le\break
\Tr(\widehat\Sigma^j(A_{\lambda}^j)^\top A_{\lambda}^j)$ while the vectors
$\mathbf{b}_{\lambda}$
are such that $A_{\lambda}^j\mathbf{b}_{\lambda}^j=(A_{\lambda
}^j)^\top \mathbf{b}_{\lambda}^j=0$. In this case, for
a temperature vector
$\bbeta=(\beta_1,\ldots,\beta_J)^\top $ and a prior $\pi$, we
define GEWA
as ${\hat{\mathbf{f}}}_\mathrm{GEWA}^j=\int_\Lambda{\tilde{\ff
}_{\lambda}}^j\hat\pi^j(d\lambda)$,
where ${\tilde{\ff}_{\lambda}}^j=(A_{\lambda}^j+(A_{\lambda
}^j)^\top -(A_{\lambda}^j)^\top A_{\lambda}^j)\YY^j+\mathbf
{b}_{\lambda}^j$
and $\hat\pi^j$ is defined by (\ref{eqdefgweights}).
Note that this setting is the grouped version of the SEWA.
\end{longlist}

%
\begin{theorem}\label{thmgroup}
Assume that $\widehat\Sigma$ is unbiased and independent of $\YY$.
Under setting~1, if $\beta_j\ge8|\!|\!|{\Sigma^j}|\!|\!|$ for all
$j=1,\ldots,J$, then
%
\begin{equation}
\label{ineqKLgroup} \quad\E \bigl[\|{\hat{\mathbf{f}}}_\mathrm {GEWA}-\ff
\|_n^2 \bigr] \leq\sum_{j=1}^J
\inf_{p_j} \biggl\{ \int_{\Lambda}\E\bigl\|\hat{
\mathbf{f}}_{\lambda}^j-\ff^j\bigr\|_n^2
p_j(d\lambda)+\frac{\beta_j}n \KL(p_j,\pi) \biggr\}.
\end{equation}
Under setting~2, this inequality holds true if $\beta_j\ge4|\!|\!
|{\Sigma^j}|\!|\!|$ for
every $j=1,\ldots,J$.
\end{theorem}

As we shall see in~Section~\ref{secmmx}, this theorem allows us to
propose an estimator of the unknown signal which is
adaptive w.r.t. the smoothness properties of the underlying signal and
achieves the minimax rates and constants over the
Sobolev ellipsoids provided that the operator $T$ is mildly ill-posed,
that is, its singular values decrease at most
polynomially.

\section{Examples of sharp oracle inequalities}\label{secexamples}

In this section we discuss consequences of the main result for specific
choices of prior measures.
For conveying the main messages of this section it is enough to focus
on settings 1 and 2 in the case of only
one group ($J=1$).

\subsection{Discrete oracle inequality}

In order to demonstrate that inequality (\ref{ineqKLgroup}) can be
reformulated in terms of an OI as defined
by (\ref{ineqoracleex}), let
us consider the case when the prior $\pi$ is discrete, that is, $\pi
(\Lambda_0)=1$ for a
countable set $\Lambda_0\subset\Lambda$, and w.l.o.g $\Lambda_0=\N
$. Then, the following result holds true.

%
\begin{proposition}\label{propdiscret}
Let $\widehat\Sigma$ be unbiased, independent of $\YY$ and $\pi$ be
supported by~$\N$. Under setting~1 with $J=1$ and
$\beta=\beta_1\ge8|\!|\!|{\Sigma}|\!|\!|$, the aggregate ${\hat
{\mathbf{f}}}_\mathrm{GEWA}$
satisfies the inequality
%
\begin{equation}
\label{ineqKLdiscrete} \E \bigl[\|{\hat{\mathbf{f}}}_\mathrm {GEWA}-\ff
\|_n^2 \bigr] \leq\inf_{\ell\in\N\dvtx\pi
_\ell>0} \biggl( \E \bigl[ \|
\hat{\mathbf{f}}_{\ell}-\ff \|_n^2 \bigr] +
\frac{\beta\log(1/\pi_\ell)}n \biggr).
\end{equation}
Furthermore, \eqref{ineqKLdiscrete} holds true under setting~2 for
$\beta\ge4|\!|\!|{\Sigma}|\!|\!|$.\vadjust{\goodbreak}
\end{proposition}
\begin{pf}
It suffices to apply Theorem~\ref{thmgroup} and to upper-bound the
right-hand side by the minimum over all Dirac
measures $p=\delta_{\ell}$ such that $\pi_\ell>0$.
\end{pf}
This inequality can be compared to Corollary 2 in~\cite{BaraudGiraudHuet10}, Section 4.3. Our result
has the advantage of having factor one in front of the expectation of
the left-hand side, while
in~\cite{BaraudGiraudHuet10} a constant much larger than 1 appears.
However, it should be noted that
the assumptions on the (estimated) noise covariance matrix are much
weaker in~\cite{BaraudGiraudHuet10}.

\subsection{Continuous oracle inequality}

It may be useful in practice to combine a family of affine estimators
indexed by an open
subset of $\R^M$ for some $M\in\N$ (\eg to build an estimator nearly as
accurate as the best kernel estimator with fixed kernel and varying
bandwidth). To state an oracle inequality
in such a ``continuous'' setup, let us denote by $d_2(\llambda
,\partial\Lambda)$ the largest real $\tau>0$ such that the ball
centered at $\llambda$ of radius $\tau$---hereafter denoted by
$B_{\llambda}(\tau)$---is included in $\Lambda$.
Let $\operatorname{Leb}(\cdot)$ be the Lebesgue measure in $\R^M$.

%
\begin{proposition}\label{prop2}
Let $\widehat\Sigma$ be unbiased, independent of $\YY$. Let $\Lambda
\subset\R^M$ be an open and bounded set and
let $\pi$ be the uniform distribution on $\Lambda$. Assume that the mapping
$\llambda\mapsto r_{\llambda}$ is Lipschitz continuous, \textit{that
is}, $|r_{\llambda'}-r_{\llambda}|\le
L_r\|\llambda'-\llambda\|_2$,
$\forall\llambda,\llambda'\in\Lambda$. Under setting~1 with $J=1$
and $\beta=\beta_1\ge8|\!|\!|{\Sigma}|\!|\!|$, the
aggregate ${\hat{\mathbf{f}}}_\mathrm{GEWA}$ satisfies the inequality
%
\begin{eqnarray}
\label{ineqKLcont} \E\|{\hat{\mathbf{f}}}_\mathrm{GEWA}-\ff
\|_n^2 &\leq& \inf_{\llambda\in\Lambda} \biggl\{ \E \bigl[\|
\hat{\mathbf {f}}_{\lambda}-\ff\|_n^2 \bigr] +
\frac{\beta M}{n} \log \biggl(\frac{\sqrt{M}}{2\min
(n^{-1},d_2(\llambda,\partial\Lambda))} \biggr) \biggr\}
\nonumber\hspace*{-30pt}
\\[-4pt]
\\[-12pt]
\nonumber
&&{}+\frac{L_r+\beta\log(\operatorname{Leb}(\Lambda))}{n}.\hspace*{-30pt}
\end{eqnarray}
Furthermore, \eqref{ineqKLcont} holds true under setting~2 for every
$\beta\ge4|\!|\!|{\Sigma}|\!|\!|$.
\end{proposition}

\begin{pf}
It suffices to apply assertion (i) of Theorem~\ref{mainthm} and to
upper-bound the right-hand side in inequality \eqref{ineqKL}
by the minimum over all measures having as density $p_{\llambda^*,\tau
^*}(\llambda)=\1_{B_{\llambda^*}(\tau^*)}(\llambda)/\operatorname
{Leb}(B_{\llambda^*}(\tau^*))$.
Choosing $\tau^*=\min(n^{-1}, d_2(\llambda^*,\partial\Lambda))$
such that $B_{\llambda^*}(\tau^*)
\subset\Lambda$, the measure $p_{\llambda^*,\tau^*}(\llambda)
\,d\llambda$ is absolutely continuous w.r.t. the uniform prior $\pi$ and the Kullback--Leibler divergence between
these two measures equals
$\log\{\operatorname{Leb}(\Lambda)/\operatorname{Leb} (B_{\llambda^*}
(\tau^*) ) \}$. Using $\operatorname{Leb} (B_{\llambda^*}
(\tau^*) )\ge
({2\tau^*}/{\sqrt{M}})^M$ and the Lipschitz condition, we get
the desired inequality.
\end{pf}

Note that it is not very stringent to require the risk function
$r_{\llambda}$ to be Lipschitz
continuous, especially since this condition needs not be satisfied
uniformly in $\ff$.
Let us consider the ridge regression: for a given design matrix $X\in
\R^{n\times p}$,
$A_{\lambda}=X(X^\top X+\gamma_n\lambda I_{n\times n})^{-1}X^\top $
and $\mathbf{b}_{\lambda}=0$ with
$\lambda\in[{\lambda_*},{\lambda^*}]$,
$\gamma_n$ being a given normalization factor typically set to $n$ or
$\sqrt{n}$, $\lambda_*>0$ and
$\lambda^*\in[\lambda_*,\infty]$. One can easily check the
Lipschitz property of the risk function
with $L_r=L_r(f)=4\lambda_*^{-1} \|\ff\|_n^2+(2/n)\Tr(\Sigma)$.

\subsection{Sparsity oracle inequality}

The continuous oracle inequality stated in the previous subsection is
well adapted to the problems
in which the dimension $M$ of $\Lambda$ is small w.r.t. the sample
size $n$ (or, more precisely, the
signal to noise ratio $n/|\!|\!|{\Sigma}|\!|\!|$). When this is not the
case, the choice of the prior should
be done more carefully. For instance, consider $\Lambda\subset\R^M$
with large $M$ under the sparsity
scenario: there is a sparse vector $\llambda^*\in\Lambda$ such that
the risk of
$\hat{\mathbf{f}}_{\llambda^*}$ is small. Then, it is natural to
choose a prior
that favors sparse
$\llambda$'s. This can be done in the same vein as in~\cite
{DalalyanTsybakov07,DalalyanTsybakov08,DalalyanTsybakov12a,DalalyanTsybakov12b},
by means of the heavy tailed prior,
%
\begin{equation}
\label{eqsparseprior} \pi(d\llambda)\propto\prod_{m=1}^M
\frac{1}{(1+|\lambda_m/\tau
|^2)^2} \1_{\Lambda}(\llambda),
\end{equation}
where $\tau>0$ is a tuning parameter.

%
\begin{proposition}\label{prop3}
Let $\widehat\Sigma$ be unbiased, independent of $\YY$. Let $\Lambda
=\R^M$ and let $\pi$ be defined by (\ref{eqsparseprior}).
Assume that the mapping $\llambda\mapsto r_{\llambda}$ is
continuously differentiable and, for some $M\times M$
matrix $\mathcal M$, satisfies
%
\begin{equation}
\label{eqLip2} r_{\llambda}-r_{\llambda'}-\nabla r_{\llambda'}^\top
\bigl(\llambda-\llambda' \bigr) \le \bigl(\llambda-
\llambda' \bigr)^\top \mathcal M \bigl(\llambda-
\llambda' \bigr)\qquad \forall\llambda,\llambda'\in\Lambda.
\end{equation}
Under setting~1 if $\beta\ge8|\!|\!|{\Sigma}|\!|\!|$, then the aggregate
${\hat{\mathbf{f}}}_\mathrm{EWA}={\hat{\mathbf{f}}}_\mathrm
{GEWA}$ satisfies
%
\begin{eqnarray}
\label{ineqKLsparse} \E \bigl[\|{\hat{\mathbf{f}}}_\mathrm {GEWA}-\ff
\|_n^2 \bigr] &\leq&\inf_{\llambda\in\R^M} \Biggl\{\E\|\hat{
\mathbf {f}}_{\lambda}-\ff\|_n^2 + \frac{4\beta}n
\sum_{m=1}^M \log \biggl(1+
\frac{|\lambda_m|}{\tau
} \biggr) \Biggr\}
\nonumber
\\[-8pt]
\\[-8pt]
\nonumber
&&{}+ \Tr(\mathcal M)\tau^2.
\end{eqnarray}
Moreover, \eqref{ineqKLsparse} holds true under setting~2 if $\beta
\ge4|\!|\!|{\Sigma}|\!|\!|$.
\end{proposition}

Let us discuss here some consequences of this sparsity oracle
inequality. First of all, consider the case of
(linearly) combining frozen estimators, that is, when
$\hat{\mathbf{f}}_{\lambda}= \sum_{j=1}^M \lambda_j \varphi_j$
with some known functions
$\varphi_j$. Then, it is clear that $r_{\llambda}-r_{\llambda
'}-\nabla r_{\llambda'}^\top (\llambda-\llambda')=2 (\llambda
-\llambda')^\top
\Phi(\llambda-\llambda')$, where $\Phi$ is the Gram matrix defined
by $\Phi_{i,j}=\langle\varphi_i|\varphi_j \rangle_n$.
So the condition in Proposition~\ref{prop3} consists in bounding the
Gram matrix of the atoms $\varphi_j$.
Let us remark that in this case---see, for instance, \cite
{DalalyanTsybakov08,DalalyanTsybakov12b}---$\Tr(\mathcal M)$ is
on the order of $M$ and the choice $\tau=\sqrt{\beta/(n M)}$ ensures
that the last term in the right-hand side of equation~(\ref{ineqKLsparse})
decreases at the parametric rate $1/n$. This is the choice we recommend
for practical applications.\vadjust{\goodbreak}

As a second example, let us consider the case of a large number of
linear estimators $\hat{\boldg}_1=G_1\YY,\ldots,\hat{\boldg
}_M=G_M\YY$
satisfying conditions of setting~1 and such that $\max_{m=1,\ldots
,M}|\!|\!|{G_m}|\!|\!|\le1$.
Assume we aim at proposing an estimator mimicking the behavior of the
best possible
convex combination of a pair of estimators chosen among $\hat{\boldg
}_1,\ldots,\hat{\boldg}_M$. This task
can be accomplished in our framework by setting $\Lambda=\R^M$ and
$\hat{\mathbf{f}}_{\lambda}=
\lambda_1 \hat{\boldg}_1+\cdots+\lambda_M\hat{\boldg}_M$, where
$\llambda=(\lambda_1,\ldots,\lambda_M)$. Remark
that if $\{\hat{\boldg}_m\}$ satisfies conditions of setting~1, so does
$\{\hat{\mathbf{f}}_{\lambda}\}$. Moreover, the mapping
$\llambda\mapsto r_{\llambda}$ is quadratic with Hessian matrix
$\nabla^2 r_{\llambda}$ given by the entries
$2\langle G_m \ff| G_{m'}\ff\rangle_n+\frac{2}{n}\Tr(G_{m'}\Sigma
G_m)$, $m,m'=1,\ldots,M$. It implies that
inequality~(\ref{eqLip2}) holds with $\mathcal M=\nabla^2 r_{\llambda}
/2$. Therefore, denoting by $\sigma_i^2$ the $i$th diagonal entry of
$\Sigma$ and setting $\ssigma=(\sigma_1,\ldots,\sigma_n)$,
we get $\Tr(\mathcal M)\le|\!|\!|{\sum_{m=1}^M G_m^2}|\!|\!| [\|\ff
\|_n^2+\|\ssigma\|_n^2 ]\le
M [\|\ff\|_n^2+\|\ssigma\|_n^2 ]$.
Applying Proposition~\ref{prop3} with $\tau=\sqrt{\beta/(n M)}$, we get
%
\begin{eqnarray}
\label{ineqKLsparse1} \E \bigl[ \|{\hat{\mathbf{f}}}_\mathrm {EWA}-\ff
\|_n^2 \bigr] &\leq&\inf_{\alpha, m,m'}  \E \bigl[\bigl\|\alpha
\hat{\boldg}_{m} +(1-\alpha)\hat{\boldg}_{m'}-\ff
\bigr\|_n^2 \bigr]
\nonumber
\\[-8pt]
\\[-8pt]
\nonumber
&&{}+ \frac{8\beta}n \log \biggl(1+ \biggl[{\frac{Mn}{\beta}}
\biggr]^{1/2} \biggr)+ \frac{\beta}{n} \bigl[\|\ff
\|_n^2+\|\ssigma\|_n^2 \bigr],
\end{eqnarray}
where the $\inf$ is taken over all $\alpha\in[0,1]$ and $m,m'\in\{
1,\ldots,M\}$.
This inequality is derived from \eqref{ineqKLsparse} by
upper-bounding the $\inf_{\llambda\in\R^M}$
by the infimum over $\llambda$'s having at most two nonzero
coefficients, $\lambda_{m^{\vphantom{1}}_0}$ and $\lambda_{m'_0}$,
that are nonnegative and sum to one: $\lambda_{m^{\vphantom
{1}}_0}+\lambda_{m'_0}=1$. To get \eqref{ineqKLsparse1}, one simply
notes that only two terms of the sum $\sum_{m}\log(1+|\lambda_{m}|\tau^{-1} )$ are nonzero and each of them is
not larger than $\log(1+\tau^{-1} )$. Thus, one can achieve
using EWA the best possible risk
over the convex combinations of a pair of linear estimators---selected
from a large (but finite) family---at
the price of a residual term that decreases at the parametric rate up
to a log factor.

\subsection{Oracle inequalities for varying-block-shrinkage
estimators}\label{sec2stein}
Let us consider now the problem of aggregation of two-block shrinkage
estimators.
This means that the constituent estimators have the following form: for
$\llambda=
(a,b,k)\in[0,1]^2\times\{1,\ldots,n\}:=\Lambda$, $\hat{\mathbf
{f}}_{\lambda}
=A_{\llambda}\YY$ where
$A_{\llambda}=\diag(a \1(i\le k)+b \1(i> k),i=1, \ldots,n
)$. Let us choose
the prior $\pi$ as uniform on $\Lambda$.

%
\begin{proposition}\label{propshrinkageoi}
Let ${\hat{\mathbf{f}}}_\mathrm{EWA}$ be the exponentially weighted
aggregate having as
constituent estimators
two-block shrinkage estimators $A_{\llambda}\YY$. If $\Sigma$ is
diagonal, then
for any $\llambda\in\Lambda$ and for any $\beta\ge8|\!|\!|{\Sigma
}|\!|\!|$,
%
\begin{equation}
\label{ineqoracleshrinkage}\qquad \E \bigl[\|{\hat{\mathbf{f}}}_\mathrm {EWA}-
\ff\|_n^2 \bigr] \leq\E \bigl[\|\hat{\mathbf{f}}_{\llambda}-
\ff\|_n^2 \bigr]+ \frac{\beta}{n} \biggl\{1+\log
\biggl(\frac{n^2\|\ff\|_n^2+n\Tr
(\Sigma)}{12\beta} \biggr) \biggr\} .
\end{equation}
\end{proposition}

In the case $\Sigma=I_{n\times n}$, this result is comparable to~\cite{Leung}, page~20, Theorem~2.49,
which\vadjust{\goodbreak} states that in the homoscedastic regression model ($\Sigma
=I_{n\times n}
$), EWA acting
on two-block positive-part James--Stein estimators satisfies, for any
$\llambda\in\Lambda$ such that
$3\le k \le n-3$ and for $\beta=8$, the oracle inequality
%
\begin{equation}
\E \bigl[\|\hat{\mathbf{f}}_{\mathrm{Leung}}-\ff\|_n^2
\bigr] \leq \E \bigl[\| \hat{\mathbf{f}}_{\llambda}-\ff\|_n^2
\bigr]+ \frac{9}{n} +\frac{8}{n}\min_{K>0} \biggl\{K\vee
\biggl( \log \frac
{n-6}{K}-1 \biggr) \biggr\}.\hspace*{-35pt}
\end{equation}

\section{Application to minimax adaptive estimation}\label{secmmx}

Pinsker proved in his celebrated paper~\cite{Pinsker80} that in the
model (\ref{eqmodel}) the minimax risk
over ellipsoids can be asymptotically attained by a linear estimator.
Let us denote by $\theta_k(\ff)=\langle\ff
|\varphi_k\rangle_n$ the coefficients of the (orthogonal) discrete
cosine\footnote{The results of this section hold true
not only for the discrete cosine transform, but also for any linear
transform $\DST$ such that $\DST\DST^\top =\DST^\top \DST
=n^{-1}I_{n\times n}$.}
(DCT) transform of $\ff$, hereafter denoted by $\DST \ff$.
Pinsker's result---restricted to Sobolev ellipsoids $\mathcal F_{\DST
}(\alpha,R)= \{\ff\in\R^n\dvtx\sum_{k=1}^n k^{2\alpha} \theta_k(\ff)^2 \le R \}$---
states that, as $n\to\infty$, the equivalences
%
\begin{eqnarray}
\inf_{\hat{\mathbf{f}}}\sup_{\ff\in\mathcal F_\DST(\alpha,R)} \E \bigl[\| \hat{\mathbf{f}}-\ff
\|_n^2 \bigr] &\sim&\inf_{A}
\sup_{f\in\mathcal F_{\DST
}(\alpha,R)} \E \bigl[\|A\YY-\ff\|_n^2 \bigr]
\\
&\sim&\inf_{w>0}\sup_{f\in\mathcal F_\DST(\alpha,R)} \E \bigl[\| A_{\alpha,w} \YY-
\ff\|_n^2 \bigr]
\end{eqnarray}
hold~\cite{Tsybakov09}, Theorem~3.2, where the first $\inf$ is
taken over all possible estimators~$\hat{\mathbf{f}}$ and
$A_{\alpha,w}=\DST^\top\diag((1-k^{\alpha}/w)_+;k=1,\ldots
,n )\DST$ is the Pinsker filter in the discrete
cosine basis. In simple words, this implies that the (asymptotically)
minimax estimator can be chosen from the quite
narrow class of linear estimators with Pinsker's filter. However, it
should be emphasized that the minimax linear
estimator depends on the parameters $\alpha$ and $R$, that are
generally unknown. An (adaptive) estimator, that
does not depend on $(\alpha,R)$ and is asymptotically minimax over a
large scale of Sobolev ellipsoids, has been
proposed by~Efromovich and Pinsker~\cite{EfromovichPinsker84}. The
next result, that is, a direct consequence
of
Theorem~\ref{mainthm}, shows that EWA with linear constituent
estimators is also asymptotically sharp adaptive
over Sobolev ellipsoids.

%
\begin{proposition}\label{propminimax}
Let $\llambda=(\alpha,w)\in\Lambda=\R_+^2$ and consider the prior
%
\begin{equation}
\label{eqpriorpinsker} \pi(d\llambda)= \frac{2n_\sigma^{-\alpha
/(2\alpha+1)}}{ (1+n_\sigma^{-\alpha
/(2\alpha+1)}w )^{3}}e^{-\alpha}\,d
\alpha\,d w,
\end{equation}
where $n_\sigma=n/\sigma^2$. Then, in model (\ref{eqmodel}) with
homoscedastic errors,
the aggregate ${\hat{\mathbf{f}}}_\mathrm{EWA}$ based on the
temperature $\beta=8\sigma^2$
and the constituent estimators
$\hat{\mathbf{f}}_{\alpha,w}=A_{\alpha,w}\YY$ (with $A_{\alpha
,w}$ being the
Pinsker filter)
is adaptive in the exact minimax sense\footnote{See~\cite{Tsybakov09}, Definition 3.8.} on the
family of classes $\{\mathcal F_\DST(\alpha,R)\dvtx\alpha>0,R>0\}$.
\end{proposition}

It is worth noting that the exact minimax adaptivity property of our
estimator ${\hat{\mathbf{f}}}_\mathrm{EWA}$ is achieved
without any tuning parameter.
All previously proposed methods that are provably adaptive in an exact
minimax sense
depend on some \mbox{parameters} such as the lengths of blocks for blockwise
Stein~\cite{CavalierTsybakov02} and
Efromovich--Pinsker~\cite{EfromovichPinsker96} estimators or the step
of discretization and the maximal value of
bandwidth~\cite{CavalierGolubevPicardTsybakov02}. Another nice
property of the estimator ${\hat{\mathbf{f}}}_\mathrm{EWA}$ is that
it does not require any pilot estimator based on the data splitting
device~\cite{GaiffasLecue11}.


We now turn to the setup of heteroscedastic regression, which
corresponds to ill-posed inverse problems as
described in Section~\ref{secill-posed}. To achieve adaptivity in the
exact minimax sense, we make use of ${\hat{\mathbf{f}}}_\mathrm{GEWA}$,
the grouped version of the exponentially weighted aggregate.\vspace*{1pt} We assume
hereafter that the matrix $\Sigma$ is diagonal
with diagonal entries $\sigma_1^2,\ldots,\sigma_n^2$ satisfying the
following property:
%
\begin{equation}
\label{cgamma} \exists\sigma_*,\gamma>0 \qquad\mbox{such that } \sigma_k^2=
\sigma_*^2 k^{2\gamma} \bigl(1+o_k(1) \bigr)\qquad
\mbox{as } k\to\infty.
\end{equation}
This kind of problems arises when $T$ is a differential operator or the
Radon transform~\cite{Cavalier08}, Section 1.3.
To handle such situations, we define the groups in the same spirit as
the weakly geometrically increasing blocks
in~\cite{CavalierTsybakov01}. Let $\nu=\nu_n$ be a positive integer
that increases as $n\to\infty$.
Set $\rho_n=\nu_n^{-1/3}$ and define
%
\begin{equation}
\label{kappaj} T_j= \cases{(1+\nu_n)^{j-1}-1,
&\quad $j=1,2,$\vspace*{2pt}
\cr
T_{j-1}+ \bigl\lfloor\nu_n
\rho_n(1+\rho_n)^{j-2} \bigr\rfloor, &\quad $j=3,4,
\ldots,$ }
\end{equation}
where $\lfloor x\rfloor$ stands for the largest integer strictly
smaller than $x$. Let $J$ be the smallest integer $j$ such that
$T_j\ge n$. We redefine $T_{J+1}=n$ and set $B_j=\{T_j+1,\ldots
,T_{j+1}\}$ for all $j=1,\ldots,J$.

%
\begin{proposition}\label{propminimax2} Let the groups $B_1,\ldots
,B_J$ be defined as above with $\nu_n$ satisfying $\log\nu_n/\log
n\to\infty$
and $\nu_n\to\infty$ as $n\to\infty$. Let $\llambda=(\alpha
,w)\in\Lambda=\R_+^2$ and consider the prior
%
\begin{equation}
\label{eqpriorpinsker2} \pi(d\llambda)= \frac{2n^{-\alpha/(2\alpha
+2\gamma+1)}}{ (1+n^{-\alpha
/(2\alpha+2\gamma+1)}w )^{3}}e^{-\alpha}\,d
\alpha\,d w.
\end{equation}
Then, in model (\ref{eqmodel}) with diagonal covariance matrix
$\Sigma=\operatorname{diag}(\sigma_k^2;1\le k\le n)$ satisfying condition
(\ref{cgamma}), the aggregate ${\hat{\mathbf{f}}}_\mathrm{GEWA}$
(under setting~1) based on
the temperatures $\beta_j=8\max_{i\in B_j}\sigma_i^2$ and the
constituent estimators\vspace*{1pt}
$\hat{\mathbf{f}}_{\alpha,w}=A_{\alpha,w}\YY$ (with $A_{\alpha
,w}$ being the
Pinsker filter) is adaptive in the exact minimax sense on the
family of classes $\{\mathcal F(\alpha,R)\dvtx\alpha>0,R>0\}$.
\end{proposition}

Note that this result provides an estimator attaining the optimal
constant in the minimax sense when the unknown signal
lies in an ellipsoid. This property holds because minimax estimators
over the ellipsoids are linear.
For other subsets of $\R^n$, such as hyper-rectangles, Besov bodies
and so on, this is not true anymore. However, as proved
by Donoho, Liu and MacGibbon~\cite{DonohoLiuMacGibbon90}, for
orthosymmetric quadratically convex sets the minimax linear estimators
have a risk which is within $25\%$ of
the minimax risk among all estimates. Therefore, following the approach
developed here, it is also possible
to prove that GEWA can lead to an adaptive estimator whose risk is
within $25\%$ of the minimax risk,
for a broad class of hyperrectangles.

%
\begin{figure}[b]

\includegraphics{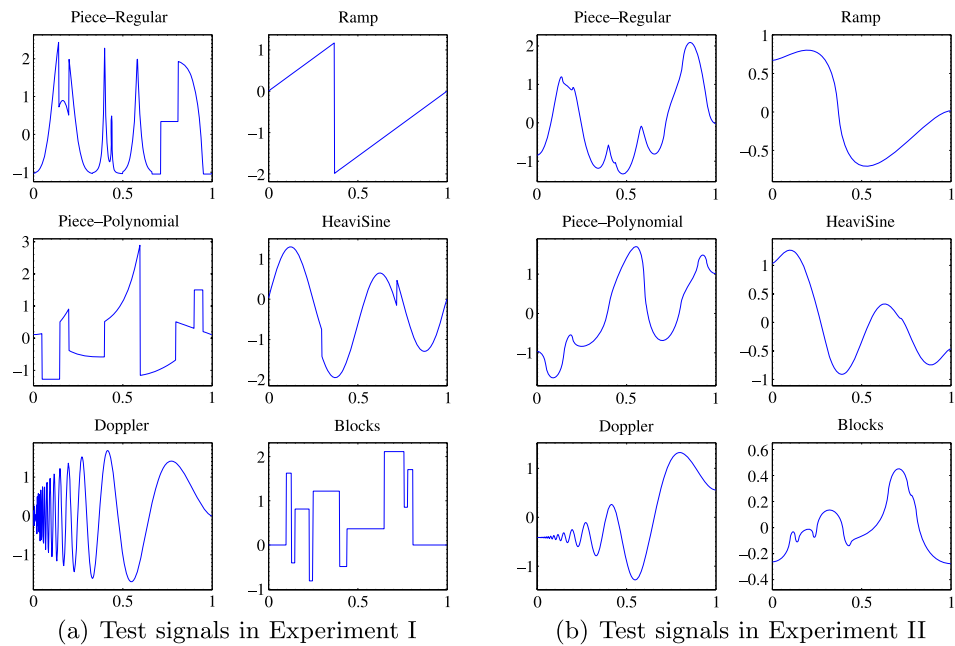}

\caption{Test signals used in our experiments: Piece-Regular, Ramp,
Piece-Polynomial, HeaviSine, Doppler and Blocks.
\textup{(a)} nonsmooth (Experiment I) and \textup{(b)} smooth (Experiment II).}\label{figtestedsignal}
\end{figure}

\section{Experiments}\label{secexperiments}

In this section we present some numerical experiments on synthetic
data, by focusing only on the
case of homoscedastic Gaussian noise ($\Sigma=\sigma^2 I_{n\times
n}$) with
known variance. A toolbox is made
available freely for download at \url{http://josephsalmon.eu/code/index\_codes.php}. Additional details and
numerical experiments can be found in \cite
{DalalyanSalmon11b,SalmonDalalyan11c}.

We evaluate different estimation routines on several 1D signals
considered as a benchmark in the literature
on signal processing~\cite{DonohoJohnstone94}. The six signals we
retained for
our experiments because of their diversity are depicted in Figure~\ref
{figtestedsignal}. Since these
signals are nonsmooth, we have also carried out experiments on their
smoothed versions obtained by
taking the antiderivative. Experiments on nonsmooth (resp., smooth)
signals are referred to as Experiment I (resp., Experiment~II).
In both cases, prior to applying
estimation routines, we normalize the (true) sampled signal to have an
empirical norm equal to one
and use the DCT denoted by
$\ttheta(\YY)= (\theta_1(\YY),\ldots,\theta_n(\YY) )^\top $.

%
%

The four tested estimation routines---including EWA---are detailed below.

\textit{Soft-Thresholding \textup{(ST)}~\cite{DonohoJohnstone94}}: For a given shrinkage parameter $t$, the
soft-thresholding
estimator is
$\widehat{\theta}_k=\sgn(\theta_k(\YY) ) (|\theta_k(\YY)|-\sigma
t )_+$.
We use the data-driven threshold minimizing the Stein unbiased risk
estimate~\cite{DonohoJohnstone95}.

\textit{Blockwise James--Stein \textup{(BJS)} shrinkage~\cite{Cai99}}: The set $\{1,\ldots,n\}$ is
partitioned into $N=[n/\log(n)]$ blocks $B_1,B_2,\ldots, B_N$ of
nearly equal size $L$. The corresponding
blocks of true coefficients $\theta_{B_k}(\ff)= (\theta_{j}(\ff
) )_{j \in B_k}$ are then
estimated by
$\widehat{\theta}_{B_k}= (1-\frac{\lambda L \sigma^2}{S_k^2(\YY)}
)_+ {\theta_{B_k}}(\YY),$ $k={1,\ldots,N}$,
with blocks of noisy coefficients ${\theta_{B_k}}(\YY)$,
$S_k^2=\|\theta_{B_k}(\YY)\|_2^2$ and $\lambda=4.50524$.

\textit{Unbiased risk estimate \textup{(URE)} minimization
with Pinsker's filters~\cite{CavalierGolubevPicardTsybakov02}}:
Pinsker filter
with data-driven parameters $\alpha$ and $w$ selected by minimizing an unbiased
estimate of the risk over a suitably chosen grid for the values of
$\alpha$ and $w$. Here, we use geometric
grids ranging from $0.1$ to $100$ for $\alpha$ and from $1$ to $n$ for $w$.

\textit{EWA on Pinsker's filters}: We consider
the same finite family of linear filters (defined by Pinsker's
filters) as in the URE routine described above. According to
Proposition~\ref{propdiscret}, this leads to an
estimator nearly as accurate as the best Pinsker's estimator in the
given family.

To report the result of our experiments, we have also computed the best
linear smoother, hereafter referred to as the oracle, based on a Pinsker
filter chosen among the candidates that we used for defining URE and
EWA. By best smoother
we mean the one minimizing the squared error (it can be computed since
we know the ground truth).
Results summarized in Table~\ref{tab1} for Experiment
I and Table~\ref{tab2} for Experiment II correspond to the average
over 1000 trials of the mean squared error
(MSE) from which we subtract the MSE of the oracle and multiply the
resulting difference by the sample size.
We report the results for $\sigma=0.33$ and for $n\in\{
2^8,2^9,2^{10},2^{11}\}$.

Simulations show that EWA and URE have very comparable performances and
are significantly more accurate than
soft-thresholding and block James--Stein (\lcf Table~\ref{tab1}) for
every size $n$ of signals considered.
Improvements are particularly important when signals have large peaks
or discontinuities. In most cases,
EWA also outperforms URE, but differences are less pronounced. One can
also observe that for smooth signals, the
difference of MSEs between EWA and the oracle, multiplied by $n$,
remains nearly constant when $n$ varies.
This is in agreement with our theoretical results in which the residual
term decreases to zero inversely
proportionally to $n$.

Of course, soft-thresholding and blockwise James--Stein procedures have
been designed for being applied
to the wavelet transform of a Besov smooth function, rather than to the
Fourier transform of a Sobolev-smooth
function. However, the point here is not to demonstrate the superiority
of EWA as compared to ST and
BJS procedures. The point is to stress the importance of having sharp
adaptivity up to an optimal constant and
not simply adaptivity in the sense of rate of convergence. Indeed, the
procedures ST and BJS are provably
rate-adaptive when applied to the Fourier transform of a Sobolev-smooth
function, but they are not sharp
adaptive---they do not attain the optimal constant---whereas EWA and
URE do attain.

%
%
%
\begin{table}
\caption{Evaluation of 4 adaptive methods on 6 (nonsmooth) signals.
For each sample size
and each method, we report the average value of $n(\MSE-\MSE_{\Oracle})$
and the corresponding standard deviation
(in parentheses), for 1000 replications of the experiment}\label{tab1}
\begin{tabular*}{\textwidth}{@{\extracolsep{\fill}}ld{2.3}d{2.3}d{2.3}d{2.3}cd{2.3}d{2.3}d{2.3}d{2.3}@{}}
\hline
\multicolumn{1}{@{}l}{$\bolds{n}$}&\multicolumn{1}{c}{\textbf{EWA}}&
\multicolumn{1}{c}{\textbf{URE}}&\multicolumn{1}{c}{\textbf{BJS}}&\multicolumn{1}{c}{\textbf{ST}}&&
\multicolumn{1}{c}{\textbf{EWA}}&\multicolumn{1}{c}{\textbf{URE}}&\multicolumn{1}{c}{\textbf{BJS}}&\multicolumn{1}{c@{}}{\textbf{ST}}\\
\hline
& \multicolumn{4}{c}{{Blocks}} & & \multicolumn{4}{c}{{Doppler}} \\
\phantom{0}256&0.051&0.245&9.617&4.846&&0.062&0.212&13.233&6.036\\
&(0.42)&(0.39)&(1.78)&(1.29)&&(0.35)&(0.31)&(2.11)&(1.23)\\
\phantom{0}512&-0.052&0.302&13.807&9.256&&-0.100&0.205&17.080&12.620\\
&(0.35)&(0.50)&(2.16)&(1.70)&&(0.30)&(0.39)&(2.29)&(1.75)\\
1024&-0.050&0.299&19.984&17.569&&-0.107&0.270&21.862&23.006\\
&(0.36)&(0.46)&(2.68)&(2.17)&&(0.35)&(0.41)&(2.92)&(2.35)\\
2048&-0.007&0.362&28.948&30.447&&-0.150&0.234&28.733&38.671\\
&(0.42)&(0.57)&(3.31)&(2.96)&&(0.34)&(0.42)&(3.19)&(3.02)\\[3pt]
& \multicolumn{4}{c}{{HeaviSine}} & & \multicolumn{4}{c}{{Piece-Regular}}
\\
\phantom{0}256&-0.060&0.247&1.155&3.966&&-0.069&0.248&8.883&4.879\\
&(0.19)&(0.42)&(0.57)&(1.12)&&(0.32)&(0.40)&(1.76)&(1.20)\\
\phantom{0}512&-0.079&0.215&2.064&5.889&&-0.105&0.237&12.147&9.793\\
&(0.19)&(0.39)&(0.86)&(1.36)&&(0.30)&(0.37)&(2.28)&(1.64)\\
1024&-0.059&0.240&3.120&8.685&&-0.092&0.291&15.207&16.798\\
&(0.23)&(0.36)&(1.20)&(1.64)&&(0.34)&(0.46)&(2.18)&(2.13)\\
2048&-0.051&0.278&4.858&12.667&&-0.059&0.283&21.543&27.387\\
&(0.25)&(0.48)&(1.42)&(2.03)&&(0.34)&(0.54)&(2.47)&(2.77)\\[3pt]
& \multicolumn{4}{c}{{Ramp}} & &
\multicolumn{4}{c}{{Piece-Polynomial}} \\
\phantom{0}256&0.038&0.294&6.933&5.644&&0.017&0.203&12.201&3.988\\
&(0.37)&(0.47)&(1.54)&(1.20)&&(0.37)&(0.37)&(1.81)&(1.19)\\
\phantom{0}512&0.010&0.293&9.712&9.977&&-0.078&0.312&17.765&9.031\\
&(0.36)&(0.51)&(1.76)&(1.67)&&(0.35)&(0.49)&(2.72)&(1.62)\\
1024&-0.002&0.300&13.656&16.790&&-0.026&0.321&23.321&17.565\\
&(0.30)&(0.45)&(2.25)&(2.06)&&(0.38)&(0.48)&(2.96)&(2.28)\\
2048&0.007&0.312&19.113&27.315&&-0.007&0.314&31.550&29.461\\
&(0.34)&(0.50)&(2.68)&(2.61)&&(0.41)&(0.49)&(3.05)&(2.95)\\
\hline
\end{tabular*}
\end{table}

%
\begin{table}
\caption{Evaluation of 4 adaptive methods on 6 smoothed signals. For
each sample size
and each method, we report the average value of $n(\MSE-\MSE_{\Oracle
})$ and the corresponding standard deviation
(in parentheses), for 1000 replications of the experiment}\label{tab2}
\begin{tabular*}{\textwidth}{@{\extracolsep{\fill}}ld{2.3}d{2.3}d{2.3}d{2.3}cd{2.3}d{2.3}d{2.3}d{2.3}@{}}
\hline
\multicolumn{1}{@{}l}{$\bolds{n}$}&\multicolumn{1}{c}{\textbf{EWA}}&
\multicolumn{1}{c}{\textbf{URE}}&\multicolumn{1}{c}{\textbf{BJS}}&\multicolumn{1}{c}{\textbf{ST}}&&
\multicolumn{1}{c}{\textbf{EWA}}&\multicolumn{1}{c}{\textbf{URE}}&\multicolumn{1}{c}{\textbf{BJS}}&\multicolumn{1}{c@{}}{\textbf{ST}}\\
\hline
& \multicolumn{4}{c}{{Blocks}} & &
\multicolumn{4}{c}{{Doppler}} \\
\phantom{0}256&0.387&0.216&0.216&2.278&&0.214&0.237&1.608&2.777\\
&(0.43)&(0.40)&(0.24)&(0.98)&&(0.23)&(0.40)&(0.73)&(1.04)\\
\phantom{0}512&0.170&0.209&0.650&3.193&&0.165&0.250&1.200&3.682\\
&(0.20)&(0.41)&(0.25)&(1.07)&&(0.20)&(0.44)&(0.48)&(1.24)\\
1024&0.162&0.226&1.282&4.507&&0.147&0.229&1.842&5.043\\
&(0.18)&(0.41)&(0.44)&(1.28)&&(0.19)&(0.45)&(0.86)&(1.43)\\
2048&0.120&0.220&1.574&6.107&&0.138&0.229&1.864&6.584\\
&(0.17)&(0.37)&(0.55)&(1.55)&&(0.20)&(0.40)&(1.07)&(1.58)\\[3pt]
& \multicolumn{4}{c}{{HeaviSine}} & &
\multicolumn{4}{c}{{Piece-Regular}}
\\
\phantom{0}256&0.217&0.207&1.399&2.496&&0.269&0.279&2.120&2.053\\
&(0.16)&(0.42)&(0.54)&(0.96)&&(0.27)&(0.49)&(1.09)&(0.95)\\
\phantom{0}512&0.206&0.221&0.024&3.045&&0.216&0.248&2.045&2.883\\
&(0.18)&(0.43)&(0.26)&(1.10)&&(0.20)&(0.45)&(1.17)&(1.13)\\
1024&0.179&0.200&0.113&3.905&&0.183&0.228&1.251&3.780\\
&(0.18)&(0.50)&(0.27)&(1.27)&&(0.20)&(0.41)&(0.70)&(1.37)\\
2048&0.162&0.189&0.421&5.019&&0.145&0.223&1.650&4.992\\
&(0.15)&(0.37)&(0.27)&(1.53)&&(0.19)&(0.42)&(1.12)&(1.42)\\[3pt]
& \multicolumn{4}{c}{{Ramp}} & &
\multicolumn{4}{c}{{Piece-Polynomial}} \\
\phantom{0}256&0.162&0.200&0.339&2.770&&0.215&0.257&1.486&2.649\\
&(0.16)&(0.38)&(0.24)&(1.00)&&(0.25)&(0.48)&(0.68)&(1.01)\\
\phantom{0}512&0.150&0.215&0.425&3.658&&0.170&0.243&1.865&3.683\\
&(0.18)&(0.38)&(0.23)&(1.20)&&(0.20)&(0.46)&(0.84)&(1.20)\\
1024&0.146&0.211&0.935&4.815&&0.179&0.236&1.547&5.017\\
&(0.18)&(0.39)&(0.33)&(1.35)&&(0.20)&(0.47)&(1.02)&(1.38)\\
2048&0.141&0.221&1.316&6.432&&0.165&0.210&2.246&6.628\\
&(0.20)&(0.43)&(0.42)&(1.54)&&(0.20)&(0.39)&(1.15)&(1.70)\\
\hline
\end{tabular*}
\end{table}

\section{Summary and future work}\label{secconclusion}

In this paper we have addressed the problem of aggregating a set of
affine estimators in the context of
regression with fixed design and heteroscedastic noise. Under some
assumptions on the constituent estimators,
we have proven that EWA with a suitably chosen temperature parameter
satisfies PAC-Bayesian type inequality,
from which different types of oracle inequalities have been deduced.
All these inequalities are with leading
constant one and rate-optimal residual term. As an application of our
results, we have shown that EWA acting
on Pinsker's estimators produces an adaptive estimator in the exact
minimax sense.

Next in our agenda is carrying out an experimental evaluation of the
proposed aggregate using the
approximation schemes described by Dalalyan and Tsybakov~\cite
{DalalyanTsybakov12b},
Rigollet and Tsybakov~\cite{RigolletTsybakov11,RigolletTsybakov11b}
and Alquier and Lounici~\cite{AlquierLounici10}, with a special focus
on the problems involving large scale data.

Although we do not assume the covariance matrix $\Sigma$ of the noise
to be known, our approach relies on an
unbiased estimator of $\Sigma$ which is independent on the observed
signal and on an upper bound on the largest
singular value of $\Sigma$. In some applications, such information may
be hard to obtain and it can be helpful to
relax the assumptions on $\widehat\Sigma$. This is another
interesting avenue for future research for which, we
believe, the approach developed by Giraud~\cite{Giraud08} can be of
valuable guidance.

\begin{appendix}\label{secproofs}
\section*{Appendix: Proofs of main theorems}
We develop now the detailed proofs of the results stated in the
manuscript.

\subsection{Stein's lemma}\label{secstein}

The proofs of our main results rely on Stein's lem\-ma~\cite{Stein73},
recalled below, providing
an unbiased risk estimate for any estimator that depends sufficiently
smoothly on the data vector $\YY$.

%
\begin{lemma}\label{lmstein}
Let $\YY$ be a random vector drawn form the Gaussian distribution
$\mathcal N_n(\ff,\Sigma)$. If the estimator
$\hat{\mathbf{f}}$ is a.e. differentiable in $\YY$ and the elements
of the
matrix $\nnabla\cdot\hat{\mathbf{f}}^\top :=(\partial_i\hat
{\mathbf{f}}_{j})$
have finite first moment, then
\[
\hat{r}= \| \YY- \hat{\ff} \|^2_n+ \frac{2}{n}\Tr
\bigl[\Sigma \bigl(\nnabla\cdot \hat{\mathbf{f}}^\top \bigr) \bigr]-
\frac{1}{n}\Tr[\Sigma]
\]
is an unbiased estimate of $r$, that is, $\E[\hat{r}]=r$.
\end{lemma}
The proof can be found in~\cite{Tsybakov09}, page 157. We apply
Stein's lemma to the affine estimators
$\hat{\mathbf{f}}_{\lambda}=A_{\lambda}\YY+\mathbf{b}_{\lambda
}$, with $A_{\lambda}$ an $n \times n$
deterministic real matrix and
$\mathbf{b}_{\lambda}\in\R^n$ a deterministic vector. We get that
if $\widehat
\Sigma$ is an unbiased estimator of $\Sigma$, then
$\hat r_\lambda^{\mathrm{unb}}= \|\YY- \hat{\mathbf
{f}}_{\lambda} \|^2_n+\frac{2}{n} \Tr[\widehat\Sigma A_{\lambda}
]-\frac{1}{n}\Tr[\widehat\Sigma]$
is an unbiased estimator of the risk
$r_\lambda=\E[\|\hat{\mathbf{f}}_{\lambda}-\ff\|_n^2]
= \|(A_{\lambda}-I_{n\times n}) \ff+\mathbf{b}_{\lambda}\|_n^2+\frac{1}{n}\Tr[A_{\lambda}\Sigma A_{\lambda}^\top ]$.

\subsection{An auxiliary result}
Prior to proceeding with the proof of the main theorems, we prove an
important auxiliary result which is the central ingredient
of the proofs for our main results.

%
\begin{lemma}\label{lmaux}
Let assumptions of Lemma~\ref{lmstein} be satisfied. Let $\{\hat
{\mathbf{f}}_{\lambda}
\dvtx\lambda\in\Lambda\}$ be a family of
estimators of $\ff$ and $\{\hat{r}_{\lambda}\dvtx\lambda\in
\Lambda\}$ a family
of risk estimates such that the mapping
$\YY\mapsto(\hat{\mathbf{f}}_{\lambda},\hat{r}_{\lambda})$ is
a.e. differentiable for every
$\lambda\in\Lambda$. Let $\hat r_\lambda^{\mathrm{unb}}$ be
the unbiased risk estimate of $\hat{\mathbf{f}}_{\lambda}$ given by
Stein's lemma.
\begin{longlist}[(1)]
\item[(1)] For every $\pi\in\mathcal{P}_\Lambda$ and for any
$\beta>0$,
the estimator ${\hat{\mathbf{f}}}_\mathrm{EWA}$ defined as the average
of $\hat{\mathbf{f}}_{\lambda}$ w.r.t. to the probability measure
\[
\hat\pi(\YY,d\lambda)=\theta(\YY,\lambda)\pi(d\lambda) \qquad\mbox {with } \theta(
\YY, \lambda)\propto{\exp \bigl\{-{n}\hat r_\lambda(\YY)/\beta \bigr\}
}
\]
admits
\[
\hat{r}_\mathrm{EWA}=\int_{\Lambda} \biggl(\hat
r_\lambda^{\mathrm
{unb}}-\|\hat{\mathbf{f}}_{\lambda}-{\hat{
\mathbf{f}}}_\mathrm {EWA}\|_n^2 -
\frac{2n}{\beta} \bigl\langle\nabla_{\YY} \hat r_\lambda|
\Sigma(\hat{\mathbf{f}}_{\lambda}- {\hat{\mathbf{f}}}_\mathrm {EWA})
\bigr\rangle_n \biggr)\hat\pi(d\lambda)
\]
as unbiased estimator of the risk.
\item[(2)] If, furthermore, $\hat r_\lambda\ge\hat r_\lambda^{\mathrm{unb}}$,
$\forall\lambda\in\Lambda$ and $\int_{\Lambda}\langle n\nabla_{\YY
} \hat r_\lambda|
\Sigma(\hat{\mathbf{f}}_{\lambda}-\break {\hat{\mathbf{f}}}_\mathrm
{EWA}) \rangle_n\hat\pi(d\lambda)\ge
- a\int_{\Lambda}\|\hat{\mathbf{f}}_{\lambda}-{\hat{\mathbf
{f}}}_\mathrm{EWA}\|_n^2\hat\pi(d\lambda)$
for some constant $a>0$, then for every $\beta\ge2a$ it holds that
%
\begin{equation}
\label{eqgen} \E \bigl[\|{\hat{\mathbf{f}}}_\mathrm{EWA}-\ff
\|_n^2 \bigr]\le \inf_{p\in\mathcal{P}_\Lambda} \biggl\{\int
_{\Lambda}\E[\hat r_\lambda] p(d\lambda)+
\frac{\beta\KL(p,\pi)}{n} \biggr\}.
\end{equation}
\end{longlist}
\end{lemma}

\begin{pf}
According to the Stein lemma, the quantity
%
\begin{equation}
\label{eqsteinhereto} \hat{r}_\mathrm{EWA}= \| \YY- {\hat{\mathbf
{f}}}_\mathrm{EWA}\|^2_n+ \frac{2}{n}\Tr
\bigl[\Sigma\bigl( \nnabla\cdot{\hat{\mathbf{f}}}_\mathrm{EWA}(\YY) \bigr)\bigr]-
\frac
{1}{n}\Tr[\Sigma]
\end{equation}
is an unbiased estimate of the risk $r_n=\E[\|{\hat{\mathbf
{f}}}_\mathrm{EWA}-\ff\|_n^2 ]$. Using simple algebra,
one checks that
%
\begin{eqnarray}
\label{eqvariancereversed} \| \YY- {\hat{\mathbf{f}}}_\mathrm {EWA}
\|^2_n&= \int_{\Lambda} \bigl( \| \YY-\hat{
\mathbf {f}}_{\lambda} \|_n^2 -\|\hat{
\mathbf{f}}_{\lambda}-{\hat{\mathbf{f}}}_\mathrm {EWA}
\|_n^2 \bigr) \hat\pi(d\lambda) .
\end{eqnarray}
By interchanging the integral and differential operators, we get the
following relation:
$\partial_{y_i} \hat{\mathbf{f}}_{ {\rm EWA},j}= \int_{\Lambda}\{
(\partial_{y_i} \hat f_{\lambda,j}(\YY) )
\theta(\YY,\lambda)+ \hat f_{\lambda,j}(\YY) (\partial_{y_i}
\theta
(\YY,\lambda) ) \} \pi(d\lambda)$.\break
Then, combining this equality with equations (\ref{eqsteinhereto}) and
(\ref{eqvariancereversed}) implies that
\begin{eqnarray*}
\hat{r}_\mathrm{EWA}=\int_{\Lambda} \bigl(\hat
r_\lambda^{\mathrm
{unb}}-\|\hat{\mathbf{f}}_{\lambda}-{\hat{
\mathbf{f}}}_\mathrm {EWA}\|_n^2 \bigr)\hat\pi(d
\lambda)+ \frac{2}{n} \int_{\Lambda}\Tr \bigl[\Sigma\hat{
\mathbf{f}}_{\lambda} \nabla_{\YY}\theta(\YY,\lambda)^\top
\bigr] \pi(d\lambda).
\end{eqnarray*}
After having interchanged differentiation and integration, we obtain
that $\int_{\Lambda}{\hat{\mathbf{f}}}_\mathrm{EWA}
(\nabla_{\YY} \theta(\YY,\lambda) )^\top \pi(d\lambda)=
{\hat{\mathbf{f}}}_\mathrm{EWA}\nabla_{\YY} (\int_{\Lambda
}\theta(\YY,\lambda) \pi
(d\lambda) )=0$ and, therefore,
we come up with the following expression for $\hat{r}_\mathrm{EWA}$:
\begin{eqnarray*}
\hat{r}_\mathrm{EWA}&=&\int_{\Lambda} \bigl( \hat
r_\lambda^{\mathrm{unb}}-\|\hat{\mathbf{f}}_{\lambda}-\hat{
\mathbf{f}}_n \|_n^2 + 2 \bigl\langle
\nabla_{\YY} \log\theta( \lambda) | \Sigma(\hat{\mathbf{f}}_{\lambda}-
{\hat{\mathbf {f}}}_\mathrm{EWA}) \bigr\rangle_n \bigr) \hat
\pi(d \lambda)
\\
&=&\int_{\Lambda} \bigl(\hat r_\lambda^{\mathrm{unb}}-\|
\hat {\mathbf{f}}_{\lambda}-{\hat{\mathbf{f}}}_\mathrm{EWA}
\|_n^2 - 2n\beta^{-1} \bigl\langle
\nabla_{\YY} \hat r_\lambda| \Sigma(\hat{\mathbf{f}}_{\lambda}-
{\hat{\mathbf{f}}}_\mathrm {EWA}) \bigr\rangle_n \bigr)\hat
\pi(d\lambda).
\end{eqnarray*}
This completes the proof of the first assertion of the lemma.

To prove the second assertion, let us observe that under the required
condition and in view of
the first assertion, for every $\beta\ge2a$ it holds that $\hat
{r}_\mathrm{EWA}
\le\int_{\Lambda}\hat r_\lambda^{\mathrm{unb}}\hat\pi(d\lambda)
\le\int_{\Lambda}\hat r_\lambda\hat\pi(d\lambda)\le\int_{\Lambda}\hat
r_\lambda\hat\pi(d\lambda)+\frac{\beta}{n}
\KL(\hat\pi,\pi)$. To conclude, it suffices to remark that $\hat
\pi$ is the probability measure minimizing the criterion
$\int_\Lambda\hat r_\lambda p(d\lambda)+\frac\beta{n}\KL(p,\pi)$
among all $p\in\mathcal{P}_\Lambda$.
Thus, for every $p\in\mathcal{P}_\Lambda$, we have
\[
\hat{r}_\mathrm{EWA}\le\int_\Lambda\hat
r_\lambda p(d\lambda)+ \frac\beta{n} \KL(p,\pi).
\]
Taking the expectation of both sides, the desired result follows.
\end{pf}

\subsection{\texorpdfstring{Proof of Theorem \protect\ref{mainthm}}{Proof of Theorem 1}}
\mbox{}

\textit{Assertion} (i).
In what follows, we use the matrix shorthand $I=I_{n\times n}$ and
${A}_\mathrm{EWA}
\triangleq\int_{\Lambda}A_{\lambda}
\hat\pi(d\lambda)$. We apply Lemma~\ref{lmaux} with $\hat
{r}_{\lambda}
=\hat r_\lambda^{\mathrm{unb}}$. To check the conditions of the
second part
of Lemma~\ref{lmaux}, note that in view of equations (\ref
{eqaffineestimatorsgal}) and (\ref{equnbiasedrisk}),
as well as the assumptions $A_{\lambda}^\top =A_{\lambda}$ and
$A_{\lambda'}\mathbf{b}_{\lambda}
=0$, we get
\begin{eqnarray*}
\nabla_{\YY} \hat r_\lambda^{\mathrm{unb}}=
\frac{2}{n}({I}- A_{\lambda})^\top ({I}-A_{\lambda})\YY-
\frac{2}{n}({I}-A_{\lambda
})^{\top} \mathbf{b}_{\lambda}=
\frac{2}{n}({I}-A_{\lambda})^2 \YY-\frac
{2}{n}
\mathbf{b}_{\lambda}.
\end{eqnarray*}
Recall now that for any pair of commuting matrices $P$ and $Q$ the identity
$(I-P)^2=(I-Q)^2+2 (I-\frac{P+Q}{2} )(Q-P)$ holds true.
Applying this identity to $P=A_{\lambda}$ and
$Q={A}_\mathrm{EWA}$ (in view of the commuting property of the
$A_{\lambda}$'s), we
get the following relation:
$ \langle(I-A_{\lambda})^2 \YY| \Sigma(A_{\lambda}- {A}_\mathrm
{EWA}) \YY
\rangle_n= \langle(I-{A}_\mathrm{EWA})^2 \YY|
\Sigma(A_{\lambda}- {A}_\mathrm{EWA}) \YY\rangle_n
-2 \langle( I-\frac{{A}_\mathrm{EWA}+A_{\lambda}}{2} ) ({A}_\mathrm
{EWA}-A_{\lambda})
\YY| \Sigma({A}_\mathrm{EWA}-A_{\lambda})\YY\rangle_n$.
When one integrates over $\Lambda$ with respect to the measure $\hat
\pi$, the term of the first scalar
product in the right-hand side of the last equation vanishes.
On the other hand,
\begin{eqnarray*}
&&\bigl\langle A_{\lambda}({A}_\mathrm{EWA}-A_{\lambda}) \YY|
\Sigma({A}_\mathrm{EWA}-A_{\lambda}) \YY \bigr\rangle_n
\\
&&\qquad= \bigl\langle A_{\lambda}({\hat{\mathbf{f}}}_\mathrm{EWA}-\hat {
\mathbf{f}}_{\lambda}) | \Sigma({\hat{\mathbf{f}}}_\mathrm {EWA}-
\hat{\mathbf{f}}_{\lambda}) \bigr\rangle_n
\\
&&\qquad= \bigl\langle({\hat{\mathbf{f}}}_\mathrm{EWA}-\hat{\mathbf
{f}}_{\lambda}) | A_{\lambda}\Sigma({\hat{\mathbf{f}}}_\mathrm
{EWA}-\hat{\mathbf{f}}_{\lambda}) \bigr\rangle_n
\\
&&\qquad= \frac1{2n}({\hat{\mathbf{f}}}_\mathrm{EWA}-\hat{\mathbf
{f}}_{\lambda})^\top (A_{\lambda}\Sigma+\Sigma
A_{\lambda}) ( {\hat{\mathbf{f}}}_\mathrm{EWA}-\hat{
\mathbf{f}}_{\lambda}) \ge0.
\end{eqnarray*}
Since positive semi-definiteness of matrices $\Sigma A_{\lambda
}+A_{\lambda}
\Sigma$ implies the one of the matrix $\Sigma{A}_\mathrm{EWA}+\Sigma
{A}_\mathrm{EWA}$,
we also have $\langle{A}_\mathrm{EWA}({A}_\mathrm{EWA}-A_{\lambda})
\YY| \Sigma({A}_\mathrm{EWA}-A_{\lambda}
)\YY\rangle_n\ge0$. Therefore,
\begin{eqnarray*}
&&\biggl\langle \biggl( I -\frac{{A}_\mathrm{EWA}+A_{\lambda}}{2} \biggr) ({A}_\mathrm{EWA}
-A_{\lambda}) \YY| \Sigma({A}_\mathrm{EWA}-A_{\lambda}) \YY
\biggr\rangle_n
\\
&&\qquad\le \bigl\langle({\hat{\mathbf{f}}}_\mathrm{EWA}-\hat{\mathbf
{f}}_{\lambda}) | \Sigma({\hat{\mathbf{f}}}_\mathrm{EWA}-\hat {
\mathbf{f}}_{\lambda}) \bigr\rangle_n
\\
&&\qquad =\bigl\|\Sigma^{1/2}({\hat{\mathbf{f}}}_\mathrm{EWA}-\hat{\mathbf
{f}}_{\lambda})\bigr\|_n^2.
\end{eqnarray*}
This inequality implies that
\begin{eqnarray*}
\int_{\Lambda} \bigl\langle n\nabla_{\YY} \hat
r_\lambda^{\mathrm
{unb}}| \Sigma(\hat{\mathbf{f}}_{\lambda}- {\hat{
\mathbf{f}}}_\mathrm{EWA}) \bigr\rangle_n\hat\pi(d\lambda )
\ge-4\int_{\Lambda}\bigl\| \Sigma^{1/2}(\hat{
\mathbf{f}}_{\lambda}-{\hat{\mathbf {f}}}_\mathrm{EWA})
\bigr\|_n^2\hat\pi(d \lambda).
\end{eqnarray*}
Therefore, the claim of Theorem~\ref{mainthm} holds true for every
$\beta\ge8|\!|\!|{\Sigma}|\!|\!|$.

\textit{Assertion} (ii).
Let now $\hat{\mathbf{f}}_{\lambda}=A_{\lambda}\YY+\mathbf
{b}_{\lambda}$ with symmetric $A_{\lambda}\preceq I_{n\times n}
$ and $\mathbf{b}_{\lambda}\in\operatorname{Ker}(A_{\lambda})$. Using
the definition
$\hat r_\lambda^{\mathrm{adj}}=\hat r_\lambda^{\mathrm{unb}}+\frac
1n\YY^\top (A_{\lambda}-A_{\lambda}^2)\YY$, one
easily checks that $\hat r_\lambda^{\mathrm{adj}}\ge\hat r_\lambda^{\mathrm{unb}}$ for every
$\lambda$ and that
\begin{eqnarray*}
\int_{\Lambda} \bigl\langle n\nabla\hat r_\lambda^{\mathrm{adj}}|
\Sigma(\hat{\mathbf{f}}_{\lambda}- {\hat{\mathbf{f}}}_\mathrm {EWA})
\bigr\rangle_n\hat\pi(d\lambda)&= &\int_{\Lambda} \bigl
\langle2(\YY -\hat{\mathbf{f}}_{\lambda}) | \Sigma(\hat{\mathbf{f}}_{\lambda}-
{\hat{\mathbf{f}}}_\mathrm {EWA}) \bigr\rangle_n\hat\pi(d
\lambda)
\\
&=&-2\int_{\Lambda}\bigl\| \Sigma^{1/2}(\hat{
\mathbf{f}}_{\lambda}- {\hat{\mathbf{f}}}_\mathrm{EWA})
\bigr\|_n^2 \hat\pi(d\lambda).
\end{eqnarray*}
Therefore, if $\beta\ge4|\!|\!|{\Sigma}|\!|\!|$, all the conditions
required in the second part of Lemma~\ref{lmaux} are fulfilled.
Applying this lemma, we get the desired result.

\subsection{\texorpdfstring{Proof of Theorem~\protect\ref{mainthmcor}}{Proof of Theorem 2}}

We apply the result of assertion (ii) of Theorem~\ref{mainthm} to
the prior
$\pi(d\lambda)$ replaced by the probability measure proportional to
$e^{(2/{\beta})\Tr[\widehat\Sigma(A_{\lambda}-A_{\lambda
}^\top A_{\lambda})]}\pi
(d\lambda)$.
This leads to
\begin{eqnarray*}
\E \bigl[\|\tilde{\mathbf{f}}_\mathrm{SEWA}-\ff
\|_n^2 \bigr] &\leq& \inf_{p \in\mathcal{P}_\Lambda} \biggl\{ \int
_{\Lambda}\E \bigl[\|\hat{\mathbf{f}}_{\lambda}-\ff
\|_n^2 \bigr]p(d\lambda)+\frac{\beta}n \KL(p,\pi)
\biggr\}
\\
& &{}+\frac{\beta}n\E \biggl[\log\int_\Lambda
e^{(2/{\beta)}\Tr[\widehat\Sigma(A_{\lambda}-A_{\lambda}^\top
A_{\lambda})]}\pi(d\lambda) \biggr].
\end{eqnarray*}
Condition (\ref{eqC}) entails that the last term is always
nonnegative and the result follows.

\subsection{\texorpdfstring{Proof of Theorem~\protect\ref{thmgroup}}{Proof of Theorem 3}}
Let us place ourselves in setting~1. It is clear that
$\E[\|{\hat{\mathbf{f}}}_\mathrm{GEWA}-\ff\|_n^2 ] = \sum_{j=1}^J
\E[\|
{\hat{\mathbf{f}}}_\mathrm{GEWA}^j-\ff^j\|_n^2 ]$.
For each $j\in\{1,\ldots,J\}$, since $\beta_j\ge8|\!|\!|{\Sigma^j}|\!|\!|$, one can apply
assertion (i) of Theorem~\ref{mainthm}, which leads to the desired
result. The case of setting~2
is handled in the same manner.
\end{appendix}

\section*{Acknowledgment}
The authors thank Pierre Alquier for fruitful discussions.

\begin{supplement}
\stitle{Proofs of some propositions}
\slink[doi]{10.1214/12-AOS1038SUPP} 
\sdatatype{.pdf}
\sfilename{aos1038\_supp.pdf}
\sdescription{In this supplement we
present the detailed proofs of Propositions~\ref{prop2}--\ref{propminimax2}.}
\end{supplement}

%

\printaddresses

\end{document}